\documentclass[letterpaper, conference]{ieeeconf}
\IEEEoverridecommandlockouts                              % This command is only
\usepackage[utf8]{inputenc}
\usepackage{textcomp}

\usepackage{graphicx}
\usepackage{amsmath}
\usepackage[version=4]{mhchem}
\usepackage{siunitx}
\usepackage{graphicx}
\usepackage{longtable,tabularx}
\usepackage[style=ieee ,sorting=none,maxbibnames=99,giveninits, doi=false]{biblatex}
\AtBeginBibliography{\scriptsize}
\newcommand{\rotation}[2]{ {\substack{#1 \\ {\boldsymbol{\longrightarrow}} \\ #2}} }
\usepackage{amsmath} % assumes amsmath package installed
\usepackage{amsfonts}
\usepackage{amsthm}

\theoremstyle{definition}

\usepackage{float} 
\usepackage[skip=0em,font=footnotesize]{caption}
\usepackage[margin = 1in]{geometry}

\usepackage{tikz}
\usetikzlibrary{shapes,arrows,calc,positioning}
\tikzstyle{bigblock} = [draw, fill=blue!20, rectangle, 
    minimum height=6em, minimum width=8em]
\tikzstyle{medblock} = [draw, fill=blue!20, rectangle, 
    minimum height=4em, minimum width=4em]    
\tikzstyle{mux} = [draw, fill=black!20, rectangle, 
    minimum height=5em, minimum width=0.1em]    
\tikzstyle{smallblock} = [draw, fill=blue!20, rectangle, 
    minimum height=2em, minimum width=3em]
    
\tikzstyle{data_block} = [draw, fill=green!20, rectangle, 
    minimum height=2em, minimum width=3em]
\tikzstyle{ops_block} = [draw, fill=blue!20, rectangle, 
    minimum height=2em, minimum width=3em]    
\tikzstyle{est_block} = [draw, fill=red!20, rectangle, 
    minimum height=2em, minimum width=3em]    
    
\tikzstyle{sum} = [draw, fill=blue!20, circle, node distance=1cm,minimum height=0.5cm]
\tikzstyle{signal} = [coordinate]
\tikzstyle{pinstyle} = [pin edge={to-,thin,black}]
\tikzstyle{block} = [draw, fill=blue!20, rectangle, 
    minimum height=3em, minimum width=9em]
\tikzstyle{blockS} = [draw, fill=blue!20, rectangle, 
    minimum height=3em, minimum width=4em]    
\tikzstyle{input} = [coordinate]
\tikzstyle{output} = [coordinate]
\usetikzlibrary{matrix}

\usetikzlibrary{positioning}
\usetikzlibrary{math}

\usepackage{hyperref}
\usepackage{xcolor}
\hypersetup{
    colorlinks,
    linkcolor={blue!100!black},
    citecolor={blue!50!black},
    urlcolor={blue!80!black}
}

\newcommand{\bc}{\begin{center}}
\newcommand{\ec}{\end{center}}
\newcommand{\benum}{\begin{enumerate}}
\newcommand{\eenum}{\end{enumerate}}
\newcommand{\nn}{\nonumber}
\newcommand{\matl}{\left[ \begin{array}}
\newcommand{\matr}{\end{array} \right]}

\renewcommand{\matl}{\begin{bmatrix}}
\renewcommand{\matr}{\end{bmatrix}}

\newcommand{\matls}{\left[ \begin{smallmatrix}}
\newcommand{\matrs}{\end{smallmatrix} \right]}
\newcommand{\isdef}{\stackrel{\triangle}{=}}

\newcommand{\vect}[1]{\overset{\rightharpoonup}{#1}}

\newcommand{\rmA}{{\rm A}}
\newcommand{\rmB}{{\rm B}}

\newcommand{\rmT}{{\rm T}}

\newcommand{\rmc}{{\rm c}}
\newcommand{\rmd}{{\rm d}}

\newcommand{\BBR}{{\mathbb R}}

\newcommand{\SB}{{\mathcal B}}

\newcommand{\SG}{{\mathcal G}}

\newcommand{\SV}{{\mathcal V}}

\newcommand{\frameddot}[2]{\stackrel{{\rm #1}\bullet \bullet}{#2}}

\newcommand{\tarrow}[1]{\overset{\rightarrow}{#1}}

\title{Adaptive Backstepping Control of a Bicopter \\ in Pure Feedback Form with Dynamic Extension}

\author{
Jhon Manuel Portella Delgado,
Mohammad Mirtaba,
and
Ankit Goel
\thanks{Jhon Manuel Portella Delgado and Mohammad Mirtaba are graduate students in the Department of Mechanical Engineering, University of Maryland, Baltimore County, 1000 Hilltop Circle, Baltimore, MD 21250. {\tt\small jportella@umbc.edu, mmirtab1@umbc.edu}}%
\thanks{Ankit Goel is an Assistant Professor in the Department of Mechanical Engineering, University of Maryland, Baltimore County,1000 Hilltop Circle, Baltimore, MD 21250. {\tt\small ankgoel@umbc.edu }}%
}

\begin{document}

\maketitle

% \begin{enumerate}
%     \item Section 3 - add Lassale statement
%     \item Section 4 - Expanbd section A and section B with commentary
%     \item Section 4 - add legends to Fig 2 Fig  3 7 and 8
%     \item Fig 10 - set ylim from 0 to 40
%     \item subscript d should be nonitalic
% \end{enumerate}

%%%%%%%%%%%%%%%%%%%%%%%%%%%%%%%%%%%%%%%%%%%%%%%%%%%%%%%%%%%%%%%%%%%%%%%%%%%%%%%%
\begin{abstract}
This paper presents a model-based, adaptive, nonlinear controller for the bicopter stabilization and trajectory-tracking problem. 
The nonlinear controller is designed using the backstepping technique. 
Due to the non-invertibility of the input map, the bicopter system is first dynamically extended.
However, the resulting dynamically extended system is in the pure feedback form with the uncertainty appearing in the input map.
The adaptive backstepping technique is then extended and applied to design the controller.
% and stabilizing the extended dynamics using the adaptive backstepping technique.
% % 
% which requires the design of additional adaptation laws.
%
The proposed controller is validated in simulation for a smooth and nonsmooth trajectory-tracking problem.

% this paper develops an adaptive nonlinear controller based on feedback linearization with a guaranteed finite time estimator for a bi-copter system. The application of IOL controllers applied to the bi-copter's dynamics suffers from invertibility problems. A solution to this problem is given by the application of time separation techniques since they treat complex systems as decoupled divisible problems. Nevertheless, they lack stability guarantees for the entire controlled closed loop. In this work, we use the dynamic extension technique to tackle this difficulty. Furthermore, since multi-copter parameters are typically time-varying or difficult to know beforehand, a finite time estimator is presented as part of an adaptive dynamic feedback linearizing controller.
% 
\end{abstract}
% \keywords{one, two, three, four}
\textit{\bf keywords:} adaptive backstepping control, dynamic extension, bicopter, multicopter.

%%%%%%%%%%%%%%%%%%%%%%%%%%%%%%%%%%%%%%%%%%%%%%%%%%%%%%%%%%%%%%%%%%%%%%%%%%%%%%%%
\section{INTRODUCTION}

Multicopters have been widely used in several engineering applications such as precision agriculture \cite{mukherjee2019}, environmental survey \cite{lucieer2014,klemas2015}, construction management \cite{li2019} and load transportation \cite{villa2020}.
This popularity has also sparked interest in broadening their capacities and applications. 
Nonetheless, due to nonlinear, time-varying, unmodeled dynamics, unknown operating environments, and evermore demanding applications, reliable multicopter control remains a challenging problem.

Various control techniques have been applied to design control systems for multicopters  \cite{nascimento2019,marshall2021,castillo2004}. Nevertheless, these techniques require prior knowledge of model parameters and, thus, are sensitive to physical model parameter uncertainty \cite{emran2018,amin2016}.
Several adaptive control techniques have been applied to address the problem of unmodeled, unknown, and uncertain dynamics, such as model reference adaptive control \cite{whitehead2010,dydek2012}, L1 adaptive control \cite{zuo2014}, adaptive sliding mode control \cite{espinoza2021trajectory,wu20221, mofid2018}, retrospective cost adaptive control \cite{goel_adaptive_pid_2021,spencer2022}.
These approaches either require an existing stabilizing controller or do not provide stability guarantees. 

Modern control architectures decompose the multicopter's nonlinear dynamics into outer-loop translational dynamics and the inner-loop rotational dynamics \cite{px4_architecture}.
Note that the translation dynamics is linear and the rotational dynamics is nonlinear. 
Stabilizing controllers are then designed for each loop separately. 
However, the cascaded multiloop architecture does not guarantee the stability of the entire closed-loop system.
The multi-loop architecture is motivated by the time separation principle, which is applicable in a scenario where each successive inner feedback loop is sufficiently faster than the previous outer loop.
A stabilizing controller can be designed for each loop with appropriate transient behavior in such a case to satisfy the time separation principle. 
% 
% These techniques are rooted in the time separation principle, which applies to the case where each successive inner feedback loop is sufficiently faster than the previous outer loop.
%
This crucial fact allows the multicopter dynamics to be decoupled and stabilizing controllers designed for each loop.
Although the controller design is considerably simplified, the closed-loop stability can not be guaranteed.

% Still and all, these cascaded controllers are arduous to tune and thus highly susceptible to failure.

This paper considers the problem of designing an adaptive controller for the fully coupled nonlinear dynamics of a multicopter. 
To focus on the controller design process, we consider a bicopter system, which is special case of a multicopter system.
A bicopter is constrained to a vertical plane and thus is modeled by a 6th-order nonlinear instead of a 12th-order nonlinear system. 
Despite the lower dimension of the state space, the 6th-order bicopter retains the complexities of the nonlinear dynamics of an unconstrained multicopter. 
% 
% 
% To remove the complexities from the typical $12th-$order nonlinear model, a $6th-$order nonlinear system is considered instead.
% 
% 
In this paper, we design an adaptive controller based on the backstepping technique \cite{krstic1995nonlinear}.
However, the classical backstepping technique can not be applied due to the input map's non invertibility.
To circumvent this problem, the bicopter dynamics is first dynamically extended \cite{descusse1985decoupling}
Although the input map of the resulting extended system is invertible, the extended dynamics can only be expressed in the pure feedback form. 
Backstepping technique has been extended to design controllers for a system in pure feedback form in \cite{zhang2017dynamic,mazenc2018backstepping,reger2019dynamic}.
An adaptive extension of the backstepping technique for a pure feedback system was presented in \cite{triska2021dynamic}.
However, \cite{triska2021dynamic} considers the case where uncertainty is in the dynamics map.
As shown in Section \ref{sec:DEBD}, the uncertainty in the bicopter dynamics appears in the input map. 
The contributions of this paper are thus
1) the design of an adaptive backstepping controller for the fully nonlinear bicopter system without decoupling the nonlinear system into simpler subsystems 
2) extension of the backstepping control of pure feedback system to the case of uncertain input maps,
and
3) validation of the proposed controller in a smooth and nonsmooth trajectory-tracking problem.

% checl line 209. please add a DE reference there ok

% ok watch

% the triska paper is yours, right?

% do you consider uncertainty in the input map in that paper?

% because that's what we have here in this paper, right?

% yes

% %
% To avoid invertibility problems of the required triangular form by the backstepping method, the $6th-$order system is first dynamically extended.
% %
% The resulting triangular structure is in pure-feedback form, where any of the intermediate virtual control or control signals are expressed in a non-affine way. This poses an extra layer of complexity due to the lack of guarantees for the existence of unique solutions. The intricacy of the problem is aggravated by considering the adaptive scenario. 
% %
% A number of methods are presented in the \textit{state of the art} to deal with the non-adaptive pure-feedback problem \cite{zhang2017dynamic,mazenc2018backstepping,reger2019dynamic}, and only a few that cope with the adaptive question such as \cite{triska2021dynamic}. 
% %
% Nevertheless, none of them has been applied to the entire closed-loop control of a fully coupled multicopter system. 
% %
% Furthermore, the adaptive scheme of this problem includes the case of virtual control unknown parameters, leveling the difficulty of the problem up.

The paper is organized as follows. 
Section \ref{sec:prob_formulation} derives the equation of motion of the bicopter system.
Section \ref{sec:adaptive_backstepping} describes the adaptive backstepping procedure to design the adaptive controller for the bicopter.
Section \ref{sec:simulations} describes a numerical simulation to validate the adaptive controller.
% the numerical simulations of the closed loop bi-copter dynamics to track a circular path and a second order Hilbert curve.
Finally, the paper concludes with a discussion of results and future research directions in section 
\ref{sec:conclusions}.

\section{Bicopter Dynamics}
\label{sec:prob_formulation}
\begin{figure}[!ht]
    \centering
\tikzset{every picture/.style={line width=0.75pt}} %set default line width to 0.75pt        

\begin{tikzpicture}[x=0.75pt,y=0.75pt,yscale=-0.8,xscale=0.8]
%uncomment if require: \path (0,300); %set diagram left start at 0, and has height of 300

%Shape: Circle [id:dp5696963057980837] 
\draw   (382.73,198.06) .. controls (384.44,196.15) and (387.38,195.98) .. (389.3,197.69) .. controls (391.22,199.41) and (391.38,202.35) .. (389.67,204.26) .. controls (387.96,206.18) and (385.02,206.35) .. (383.1,204.63) .. controls (381.18,202.92) and (381.02,199.98) .. (382.73,198.06) -- cycle ;
%Straight Lines [id:da8445085480646768] 
\draw    (385.37,137.84) -- (385.88,195.97) ;
\draw [shift={(385.35,135.84)}, rotate = 89.5] [fill={rgb, 255:red, 0; green, 0; blue, 0 }  ][line width=0.08]  [draw opacity=0] (12,-3) -- (0,0) -- (12,3) -- cycle    ;
%Straight Lines [id:da5402004060683299] 
\draw    (449.36,199.25) -- (391.23,200.02) ;
\draw [shift={(451.36,199.23)}, rotate = 179.24] [fill={rgb, 255:red, 0; green, 0; blue, 0 }  ][line width=0.08]  [draw opacity=0] (12,-3) -- (0,0) -- (12,3) -- cycle    ;
%Shape: Circle [id:dp5747098672764621] 
\draw  [fill={rgb, 255:red, 0; green, 0; blue, 0 }  ,fill opacity=1 ] (384.93,201.16) .. controls (384.93,200.46) and (385.5,199.9) .. (386.2,199.9) .. controls (386.9,199.9) and (387.47,200.46) .. (387.47,201.16) .. controls (387.47,201.86) and (386.9,202.43) .. (386.2,202.43) .. controls (385.5,202.43) and (384.93,201.86) .. (384.93,201.16) -- cycle ;
%Straight Lines [id:da7824424658656131] 
\draw    (322.35,83.75) -- (220.35,195.75) ;
\draw [shift={(271.35,139.75)}, rotate = 132.32] [color={rgb, 255:red, 0; green, 0; blue, 0 }  ][fill={rgb, 255:red, 0; green, 0; blue, 0 }  ][line width=0.75]      (0, 0) circle [x radius= 1.34, y radius= 1.34]   ;
%Straight Lines [id:da6481972501845306] 
\draw    (200.75,177.75) -- (220.35,195.75) ;
%Shape: Ellipse [id:dp22032720421843632] 
\draw   (183.73,197.05) .. controls (181.5,194.8) and (183.84,188.88) .. (188.95,183.81) .. controls (194.05,178.75) and (200,176.46) .. (202.22,178.71) .. controls (204.45,180.95) and (202.11,186.88) .. (197.01,191.94) .. controls (191.9,197.01) and (185.95,199.29) .. (183.73,197.05) -- cycle ;
%Shape: Ellipse [id:dp42306385916644906] 
\draw   (202.22,178.71) .. controls (200,176.46) and (202.33,170.54) .. (207.44,165.47) .. controls (212.55,160.41) and (218.49,158.12) .. (220.72,160.37) .. controls (222.94,162.61) and (220.61,168.54) .. (215.5,173.6) .. controls (210.39,178.66) and (204.45,180.95) .. (202.22,178.71) -- cycle ;
%Straight Lines [id:da2693350638942833] 
\draw    (302.75,65.75) -- (322.35,83.75) ;
%Shape: Ellipse [id:dp5012520329793704] 
\draw   (285.73,85.05) .. controls (283.5,82.8) and (285.84,76.88) .. (290.95,71.81) .. controls (296.05,66.75) and (302,64.46) .. (304.22,66.71) .. controls (306.45,68.95) and (304.11,74.88) .. (299.01,79.94) .. controls (293.9,85.01) and (287.95,87.29) .. (285.73,85.05) -- cycle ;
%Shape: Ellipse [id:dp1197327198597613] 
\draw   (304.22,66.71) .. controls (302,64.46) and (304.33,58.54) .. (309.44,53.47) .. controls (314.55,48.41) and (320.49,46.12) .. (322.72,48.37) .. controls (324.94,50.61) and (322.61,56.54) .. (317.5,61.6) .. controls (312.39,66.66) and (306.45,68.95) .. (304.22,66.71) -- cycle ;
%Straight Lines [id:da9675229957578317] 
\draw  [dash pattern={on 4.5pt off 4.5pt}]  (331.49,138.96) -- (271.35,139.75) ;
%Shape: Circle [id:dp6140250358557469] 
\draw  [fill={rgb, 255:red, 0; green, 0; blue, 0 }  ,fill opacity=1 ] (326.07,201.55) .. controls (326.07,200.98) and (326.53,200.51) .. (327.1,200.51) .. controls (327.68,200.51) and (328.14,200.98) .. (328.14,201.55) .. controls (328.14,202.12) and (327.68,202.59) .. (327.1,202.59) .. controls (326.53,202.59) and (326.07,202.12) .. (326.07,201.55) -- cycle ;
%Straight Lines [id:da3137247335444493] 
\draw [line width=1.5]    (410.42,70.59) -- (410.12,109.92) ;
\draw [shift={(410.09,113.92)}, rotate = 270.44] [fill={rgb, 255:red, 0; green, 0; blue, 0 }  ][line width=0.08]  [draw opacity=0] (8.75,-4.2) -- (0,0) -- (8.75,4.2) -- (5.81,0) -- cycle    ;
%Straight Lines [id:da6229937781540653] 
\draw [line width=1.5]    (196.22,172.04) -- (172.07,149.28) ;
\draw [shift={(169.16,146.54)}, rotate = 43.29] [fill={rgb, 255:red, 0; green, 0; blue, 0 }  ][line width=0.08]  [draw opacity=0] (8.75,-4.2) -- (0,0) -- (8.75,4.2) -- (5.81,0) -- cycle    ;
%Straight Lines [id:da12253033078799613] 
\draw [line width=1.5]    (298.42,60.42) -- (274.27,37.66) ;
\draw [shift={(271.35,34.92)}, rotate = 43.29] [fill={rgb, 255:red, 0; green, 0; blue, 0 }  ][line width=0.08]  [draw opacity=0] (8.75,-4.2) -- (0,0) -- (8.75,4.2) -- (5.81,0) -- cycle    ;
%Shape: Circle [id:dp662955368321148] 
\draw   (192.75,92.55) .. controls (192.54,89.99) and (194.46,87.75) .. (197.02,87.55) .. controls (199.58,87.35) and (201.82,89.26) .. (202.02,91.83) .. controls (202.22,94.39) and (200.31,96.63) .. (197.75,96.83) .. controls (195.18,97.03) and (192.95,95.11) .. (192.75,92.55) -- cycle ;
%Straight Lines [id:da15169212283243372] 
\draw    (151.07,49) -- (193.41,88.83) ;
\draw [shift={(149.61,47.63)}, rotate = 43.25] [fill={rgb, 255:red, 0; green, 0; blue, 0 }  ][line width=0.08]  [draw opacity=0] (12,-3) -- (0,0) -- (12,3) -- cycle    ;
%Straight Lines [id:da36046097574097] 
\draw    (239.68,45.24) -- (200.04,87.77) ;
\draw [shift={(241.05,43.78)}, rotate = 132.99] [fill={rgb, 255:red, 0; green, 0; blue, 0 }  ][line width=0.08]  [draw opacity=0] (12,-3) -- (0,0) -- (12,3) -- cycle    ;
%Shape: Circle [id:dp04796725084777309] 
\draw  [fill={rgb, 255:red, 0; green, 0; blue, 0 }  ,fill opacity=1 ] (196.51,93.1) .. controls (196,92.62) and (195.99,91.82) .. (196.47,91.31) .. controls (196.95,90.81) and (197.75,90.79) .. (198.26,91.27) .. controls (198.77,91.76) and (198.78,92.56) .. (198.3,93.07) .. controls (197.82,93.57) and (197.01,93.59) .. (196.51,93.1) -- cycle ;
%Shape: Arc [id:dp9925831036411028] 
\draw  [draw opacity=0] (295.51,113.84) .. controls (295.56,113.84) and (295.62,113.84) .. (295.67,113.84) .. controls (301.93,113.99) and (306.85,120.97) .. (306.66,129.45) .. controls (306.57,133.33) and (305.43,136.85) .. (303.61,139.52) -- (295.32,129.19) -- cycle ; \draw    (295.67,113.84) .. controls (301.93,113.99) and (306.85,120.97) .. (306.66,129.45) .. controls (306.59,132.56) and (305.85,135.43) .. (304.62,137.82) ; \draw [shift={(303.61,139.52)}, rotate = 291.91] [fill={rgb, 255:red, 0; green, 0; blue, 0 }  ][line width=0.08]  [draw opacity=0] (7.2,-1.8) -- (0,0) -- (7.2,1.8) -- cycle    ; \draw [shift={(295.51,113.84)}, rotate = 24.68] [fill={rgb, 255:red, 0; green, 0; blue, 0 }  ][line width=0.08]  [draw opacity=0] (7.2,-1.8) -- (0,0) -- (7.2,1.8) -- cycle    ;

% Text Node
\draw (413.67,79.33) node [anchor=north west][inner sep=0.75pt]   [align=left] {{\small $\vect g$}};
% Text Node
\draw (453,190) node [anchor=north west][inner sep=0.75pt]   [align=left] {{\small $\hat \imath _{\rm{A}}$}};
% Text Node
\draw (375,120) node [anchor=north west][inner sep=0.75pt]   [align=left] {{\small $\hat \jmath _{\rm{A}}$}};
% Text Node
\draw (234,31) node [anchor=north west][inner sep=0.75pt]  [rotate=-313.75] [align=left] {{\small $\hat \imath _{\rm B}$}};
% Text Node
\draw (129.9,40) node [anchor=north west][inner sep=0.75pt]  [rotate=-313.75] [align=left] {{\small $\hat \jmath _{\rm B}$}};
% Text Node
\draw (287,20.33) node [anchor=north west][inner sep=0.75pt]   [align=left] {{\small ${\vect f}_2$}};
% Text Node
\draw (184,135.33) node [anchor=north west][inner sep=0.75pt]   [align=left] {{\small ${\vect f}_1$}};
% Text Node
\draw (260.67,125) node [anchor=north west][inner sep=0.75pt]   [align=left] {{\small c}};
% Text Node
\draw (327.33,200) node [anchor=north west][inner sep=0.75pt]   [align=left] {{\small \textit{w}}};
% Text Node
\draw (308.33,115) node [anchor=north west][inner sep=0.75pt]   [align=left] {{\small $\theta$}};
\end{tikzpicture}
% \vspace*{0.5mm}
\caption{Bicopter configuration considered in this paper. The bicopter is constrained to the $\hat \hat \imath _\rmA-\hat \hat \jmath _\rmA$ plane and rotates about the $\hat k_\rmA$ axis of the inertial frame $\rm F_A.$  }
    \label{fig:Bicopter}
    \vspace{-2em}
\end{figure}
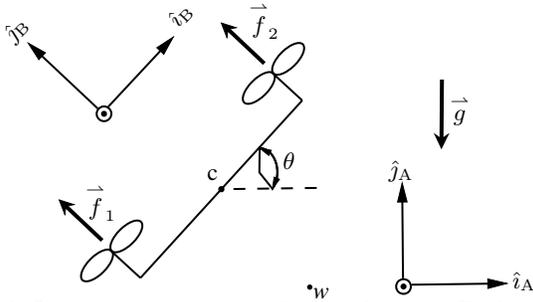

Let ${\rm F_A} = \{ \hat \imath_\rmA, \hat \jmath_\rmA, \hat k_\rmA\}$ be an inertial frame and let ${\rm F_B} = \{ \hat \imath_\rmB, \hat \jmath_\rmB, \hat k_\rmB\}$ be a frame fixed to the bicopter $\SB$ as shown in Figure \ref{fig:Bicopter}.
The bicopter $\SB$ is constrained to move in the $\hat \imath_\rmA - \hat \jmath_\rmA$ plane. 
Note that $\rm F_B$ is obtained by rotating it about the $\hat k_\rmA$ axis of $\rm F_A$ by $\theta,$ and thus 
\begin{align}
    {\rm F_A}
        \rotation{\theta}{3}
    {\rm F_B}.
\end{align}

\hfill\\

Letting $\rmc$ denote the center of mass of the bicopter and $w$ denote a fixed point on Earth, it follows from Newton's second law that 
\begin{align}
    m
    \frameddot{A}{\vect r}_{\rmc/w} 
        &=
            m \vect g + \vect f,
        \label{eq:N2L_pos}
\end{align}
where 
$m$ is the mass of the bicopter,
$\vect g$ is the acceleration due to gravity, 
and $\vect f$ is the total force applied by the propellers to the bicopter. 
Letting $\vect f_1 = f_1 \hat \jmath_\rmB$ and $\vect f_2 = f_2 \hat \jmath_\rmB$ denote the forces applied by the two propellers, it follows that $\vect f =  f_1 \hat j_{\rm B} + f_2 \hat j_\rmB $. Writing
$\vect r_{\rmc/w} 
        =
            r_1 \hat \imath_\rmA +
            r_2 \hat \jmath_\rmA$ yields
\begin{align}
    m \ddot r_1 &= -(f_1 + f_2) \sin \theta, 
    \label{eq:eom_r1}
    \\
    m \ddot r_2 &= (f_1 + f_2) \cos \theta - m g.
    \label{eq:eom_r2}
\end{align}
Next, it follows from Euler's equation that 
\begin{align}
    \tarrow J_{\SB/\rmc} \frameddot{A}{\vect \omega}_{\rm B/A} 
        =
            \vect M_{\SB/\rmc},
    \label{eq:EulersEqn}
\end{align}
where $\tarrow J_{\SB/\rmc}$ is the physical inertia matrix and  
$\vect M_{\SB/\rmc}$ is the moment applied to $\SB$ about the point $\rmc. $
Note that  $\tarrow J_{\SB/\rmc} \frameddot{A}{\vect \omega}_{\rm B/A}  = J \ddot \theta \hat k_\rmB $
and $\vect M_{\SB/\rmc} = \ell (f_2 - f_1) \hat k_\rmB,$ where $\ell$ is the length of the bicopter arm, 
and thus it follows from \eqref{eq:EulersEqn} that
\begin{align}
    J \ddot \theta = \ell (f_2 - f_1). 
    \label{eq:eom_theta}
\end{align}

Defining
\begin{align}
    F &\isdef f_1 + f_2, \label{control_signals_1}\\ 
    M & \isdef (-f_1+f_2)\ell,
    \label{control_signals_2}
\end{align}
it follows that 
\begin{align}
    m \ddot r_1 &= -F \sin \theta, 
    \label{eq:eom_r1_2}
    \\
    m \ddot r_2 &= F \cos \theta - m g,
    \label{eq:eom_r2_2}
    \\
    J \ddot \theta &= M.
    \label{eq:eom_theta_2}
\end{align}

\section{Adaptive Backstepping Control}
\label{sec:adaptive_backstepping}
In this section, we construct an adaptive backstepping controller to stabilize the bicopter system. 
To do so, the equations of motion, given by \eqref{eq:eom_r1_2}-\eqref{eq:eom_theta_2}, must be reformulated into the strict feedback form to apply the classical backstepping. 
However, as shown below, the classical backstepping technique is not applicable due to the singularity of the input map. 
% 

% Next, we write the equations of motion of the bicopter, given by \eqref{eq:eom_r1_2}-\eqref{eq:eom_theta_2}, in the triangular state-space form to develop a backstepping-based controller. 
% % 
% However, as shown below, the classical backstepping technique is not applicable due to the singularity of the input map. 

\subsection{Strict Feedback Form of Bicopter Dynamics}

Defining
\begin{align}
    x_1 \isdef \matl
                 r_1\\
                 r_2\\
                 \theta
               \matr, 
    x_2 \isdef \matl
                 \dot r_1\\
                 \dot r_2\\
                 \dot \theta
               \matr,
    u
        &\isdef \matl
                    F \\
                    M
                 \matr,
    % \label{control_signals}
    \Theta
        \isdef \matl
                    m^{-1}\\
                    J^{-1}
                  \matr\label{linear_parameters}
\end{align}
it follows that 
\begin{align}
    \dot x_1 &= x_2, 
    \\
    \dot x_2
        &= 
            f_2(x_1,x_2) + g_2(x_1,x_2) {\rm  diag}(\Theta)\,u,
\end{align}
where 
\begin{align}
    f_2(x_1,x_2) 
        \isdef 
            \matl 
                0\\
                -g\\ 
                0
            \matr
    ,\quad 
    g_2(x_1,x_2)  
        \isdef 
            \matl 
                -\sin(x_{1,3}) & 0\\
                \cos(x_{1,3}) & 0\\
                0 & 1
            \matr.
\end{align}
Note that the map $g_2$ is not invertible; hence, classical backstepping \cite[p.~29]{krstic1995nonlinear} can not be applied to this system. 

% does not fulfill assumption \ref{triangular_system_assumption} due to the non-invertibility of $g_2(x_1,x_2)$.

\subsection{Dynamic Extension of Bicopter Dynamics}
\label{sec:DEBD}
To circumvent the problem of the singularity of the input map, we extend the dynamics as shown below. 
Redefining
\begin{align}
    x_1 \isdef \matl
                 r_1\\
                 r_2
               \matr, \
    x_2 \isdef \matl
                 \dot r_1\\
                 \dot r_2
               \matr, \
    x_3 \isdef \matl
                 F\\
                 \theta
               \matr, \
    x_4 \isdef \matl
                 \dot F\\
                 \dot \theta
               \matr,
\end{align}
and 
\begin{align}
    u
        \isdef 
            \matl 
                \ddot F \\
                M
            \matr,\quad
    \Theta
        \isdef
            \matl
                \Theta_1\\
                \Theta_2
            \matr
            =
            \matl
                m^{-1}\\
                J^{-1}
            \matr,
\end{align}
% Taking into account definitions \eqref{linear_parameters} and \eqref{control_signals_dynamic_extnsion}, and
it follows from the equations of motion \eqref{eq:eom_r1}, \eqref{eq:eom_r2}, and \eqref{eq:eom_theta} that
% \begin{align}
%     \dot x_1 &= x_2\\
%     \dot x_2 &= \matl
%                     0\\
%                     -g
%                 \matr + \matl
%                           -\sin(\theta)F\\
%                           \cos(\theta)F
%                         \matr
%                         \Theta_1\\
%     \dot x_3 &= x_4\\
%     \dot x_4 
%         &= 
%             % \matl
%             %     1 & 0\\
%             %     0 & 1
%             % \matr
%             \matl
%                         1 & 0\\
%                         0 & \Theta_2
%                     \matr\,u.
% \end{align}
% 
% which can be written as
\begin{align}
    \dot x_1 &= x_2,
    \label{eq:x1dot}
    \\
    \dot x_2 &= f_2(x_1, x_2) + g_2(x_3) \Theta_1 ,
    \label{eq:x2dot}
    \\
    \dot x_3 &= x_4,
    \label{eq:x3dot}
    \\
    \dot x_4
        &=
            % g_4(x) 
            {\rm diag} (1, \Theta_2) u \label{eq:x4dot},
\end{align}
where 
\begin{align}
    f_2(x_1, x_2) 
        &\isdef 
            \matl
                0\\
                -g
            \matr,
            % \\
            \quad 
    g_2(x_3) 
        \isdef
            \matl
                -\sin(x_{3,2})x_{3,1}\\
                \cos(x_{3,2}) x_{3,1}
            \matr.
\end{align}
% Note that the system \eqref{eq:x1dot}-\eqref{eq:x4dot} is in pure feedback form since $x_3$ appears non-affinely in \eqref{eq:x2dot}. 
Note that the dynamic extension of the bicopter system, given by \eqref{eq:x1dot}-\eqref{eq:x4dot}, is not in the desired strict feedback form. 
In fact, the extended dynamics is in pure feedback form since $x_3$ appears non-affinely in \eqref{eq:x2dot}. 
% , the bicopter dynamics is no longer in the desirable strict feedback form. 
However, since the input map ${\rm diag} (1, \Theta_2)$ is invertible, and thus the backstepping technique can now be applied to the extended system to construct a controller.
% However, t

\subsection{Adaptive Backstepping-based Control of Bicopter}
\label{sec:Aadaptive backstepping}
Next, we construct an adaptive controller based on the backstepping technique.
This approach is based on the process described in \cite[p.~29]{krstic1995nonlinear}.

\subsubsection{$e_1$ Stabilization}
Let $x_{1\rmd} = \matl r_{\rmd1} & r_{\rmd2} \matr^\rmT$ denote the desired trajectory, and define the tracking error $e_1 \isdef x_1 - x_{1\rmd}.$ 
Consider the function
\begin{align}
    V_1 
    % (x_1)
        \isdef  
            \frac{1}{2} e_1^\rmT e_1.
    \label{eq:V1}
\end{align}
Differentiating \eqref{eq:V1} and using \eqref{eq:x1dot} yields
\begin{align}
    \dot V_1
        =
            e_1^\rmT \dot x_1
        =
            e_1^\rmT x_2 .
\end{align}

Note that if $x_2 = -k_1 e_1,$ where $k_1>0,$ then $\dot V_1 < 0.$    
However, $x_2$ is not the control signal and thus can not be chosen arbitrarily.
Instead, following the backstepping process, we design a control law that yields the desired $x_2$ response, as shown below.

\subsubsection{$e_2$ Stabilization}

Next, consider the function
\begin{align}
    V_2 
    % (x_1, x_2)
        \isdef 
            V_1 + 
            \frac{1}{2} 
            e_2^\rmT 
            e_2,
    \label{eq:V2}
\end{align}
where 
\begin{align}
    e_2 &\isdef x_2 - x_{2\rmd}, \label{eq:e2}\\
    x_{2\rmd} &\isdef -k_1 e_1.
    \label{eq:x2d}
\end{align}
Differentiating \eqref{eq:V2} and using \eqref{eq:x2dot} yields
\begin{align}
    \dot V_2
        &=
        %     \dot V_1 +  
        %     e_2^\rmT \dot e_2
        % \nn \\
        % &=
        %     x_1^\rmT x_2 +
        %     e_2^\rmT ( f_2 + g_2 \Theta_1 + k_1 \dot x_1)
        % \nn \\
        % &=
        %     x_1^\rmT x_2 +
        %     e_2^\rmT ( f_2 + g_2 \Theta_1 + k_1 x_2)
        % \nn \\
        % &=
        %     x_1^\rmT x_{2\rmd}
        %     +
        %     x_1^\rmT e_2         
        %     +
        %     e_2^\rmT ( f_2 + g_2 \Theta_1 + k_1 x_2)
        % \nn \\
        % &=
            -k_1 e_1^\rmT e_1
            +
            e_2^\rmT 
            % (x_1 + f_2 + g_2 \Theta_1 \nn\\
            % &\quad + k_1 x_2)
            \xi_2,
    \label{eq:V2dot}
\end{align}
where 
\begin{align}
    \xi_2 
        \isdef 
            x_1 + f_2 + g_2 \Theta_1 + k_1 x_2.
\end{align}
% Note that if $x_2 = -k_1 x_1,$ then it follows that $\dot V_1 < 0.$

Note that if 
\begin{align}
    g_2(x_3) \Theta_1 = -e_2 - e_1 - f_2 - k_1 x_2 ,
\end{align}
then $\dot V_2 < 0.$

Let $\hat \Theta_1$ be an estimate of $\Theta_1,$ and define 
$\tilde \Theta_1 \isdef \hat \Theta_1  - \Theta_1.$
Next, define
\begin{align}
    x_{3\rmd} 
        &\isdef
        %     x_3 - (x_1 + f_2 + g_2 \hat \Theta_1 + k_1 x_2)
        %     - k_2e_2
        % \nn \\
        % &=
            x_3 - (g_2 \tilde \Theta_1  +
            % x_1 + f_2  + g_2 \Theta_1  + k_1 x_2
            \xi_2 
            )
            % \nn \\ &\quad
            - k_2e_2,
            \label{eq:x3d}
    \\
    e_3 &\isdef x_3 - x_{3\rmd},
\end{align}
where $k_2>0,$
which implies that 
\begin{align}
    % x_1 + f_2  + g_2 \Theta_1  + k_1 x_2
    \xi_2 
        &=
            e_3 - g_2 \tilde \Theta_1 
            % \nn \\ &\quad
            - k_2e_2
        \in \BBR^2.
    \label{eq:V2dot_RHS_term}
\end{align}

Substituting \eqref{eq:V2dot_RHS_term} in \eqref{eq:V2dot} yields
\begin{align}
    \dot V_2
        &=
            -k_1 e_1^\rmT e_1
            - k_2 e_2^\rmT e_2
            % \nn \\ &\quad
            +
            e_2^\rmT 
            \left( e_3 + g_2 ( \Theta_1  - \hat \Theta_1) \right).        
    \label{eq:V2dot_2}
\end{align}

\subsubsection{$\Theta_1$ Adaptation}
Next, consider the function
\begin{align}
    \SV_{2} 
    % (x_1, x_2, \hat \Theta_1)
        \isdef 
            V_2 + 
            \frac{1}{2} \gamma_1^{-1} (\Theta_1  - \hat \Theta_1) ^\rmT ( \Theta_1  - \hat \Theta_1)
    \label{eq:SV2}
\end{align}
where $\gamma_1 > 0.$
Differentiating \eqref{eq:SV2} yields
\begin{align}
    \dot \SV_{2} 
        % &= 
        %     -k_1 e_1^\rmT e_1
        %     - k_2 e_2^\rmT e_2
        %     \nn \\ &\quad
        %     +
        %     e_2^\rmT 
        %     \left( e_3 + g_2 ( \Theta_1  - \hat \Theta_1) \right)
        %     \nn \\ &\quad
        %     -
        %     \gamma_1^{-1} (\Theta_1  - \hat \Theta_1) ^\rmT \dot {\hat \Theta}_1 
        % \nn \\
        &=
            -k_1 e_1^\rmT e_1
            - k_2 e_2^\rmT e_2
            % \nn \\ &\quad
            +
            e_2^\rmT  e_3 +
            \nn \\ &\quad
            ( \Theta_1  - \hat \Theta_1) ^\rmT
            \left( g_2^\rmT e_2
            -
            \gamma_1^{-1}  \dot {\hat \Theta}_1 \right).
    % \label{eq:SV2}
\end{align}
Letting
\begin{align}
    \dot {\hat \Theta}_1
        =
            \gamma_1 g_2^\rmT e_2
    \label{eq:hat_Theta_1_dot}
\end{align}
yields
\begin{align}
    \dot \SV_{2} 
        &=
            -k_1 e_1^\rmT e_1
            - k_2 e_2^\rmT e_2
            % \nn \\ &\quad
            +
            e_2^\rmT  e_3 .
    % \label{eq:SV2}
\end{align}

\subsubsection{$e_3$ Stabilization}

Next, consider the function
\begin{align}
    V_3
        \isdef
            \SV_2 
            + 
            \frac{1}{2} e_3^\rmT e_3.
    \label{eq:V3}
\end{align}
Note that
\begin{align}
    \dot\xi_2 &= x_2 + \SG_2 \Theta_1 + k_1(f_2 + g_2\Theta_1),
    \label{eq:alpha1_dot}
\end{align}
where 
\begin{align}
    \SG_2 \isdef {\partial_{x_3} g_2}
            =
            \matl
                                                        -\cos(x_{3,2})x_{3,1} & -\sin(x_{3,2})\\
                                                        -\sin(x_{3,2})x_{3,1} & \cos(x_{3,2})
                                                      \matr.
\end{align}
Differentiating \eqref{eq:V3} and using \eqref{eq:x3dot}, \eqref{eq:x3d} and \eqref{eq:alpha1_dot} yields
\begin{align}
    \dot V_3
%         &=
%             \dot \SV_{2} + e_3^\rmT (\dot {x}_3 - \dot{x}_{3d})
% %     \label{eq:V3dot}
% % \end{align}
% % Using \eqref{eq:x3dot} and the derivative of \eqref{eq:x3d} in \eqref{eq:V3dot} gives in
% % \begin{align}
% %     \dot V_3
%         \nn \\
%         &=
%            -k_1 e_1^\rmT e_1
%             - k_2 e_2^\rmT e_2
%             \nn \\ &\quad
%             + e_3^\rmT \Big(x_2-x_{2\rmd} + x_4 - x_4 + x_2 \nn\\
%             & \quad + \SG_2 \dot x_3\hat{\Theta}_1 + g_2\dot{\hat{\Theta}}_1 + k_1(f_2+g_2\Theta_1) \nn\\
%             & \quad - k_2(f_2+g_2\Theta_1) + k_1k_2x_2\Big) 
%         \nn\\
%         &=
%             -k_1 e_1^\rmT e_1
%             - k_2 e_2^\rmT e_2
%             \nn \\ &\quad
%             + e_3^\rmT \Big(2x_2 + k_1 e_{1} + \SG_2 \dot x_3\hat{\Theta}_1 \nn\\
%             & \quad + g_2\dot{\hat{\Theta}}_1 + \kappa_{12}(f_2+g_2\Theta_1) + k_1k_2x_2\Big)
%         \nn\\
        % &=
        %     -k_1 e_1^\rmT e_1
        %     - k_2 e_2^\rmT e_2
        %     \nn \\ &\quad
        %     + e_3^\rmT \Big(2x_2 + k_1 e_{1} + \SG_2 \dot x_3\hat{\Theta}_1 \nn\\
        %     & \quad + g_2\dot{\hat{\Theta}}_1 + \kappa_{12}
        %     (f_2+g_2 \hat \vartheta_1 - g_2\hat \vartheta_1 + g_2\Theta_1) + k_1k_2x_2\Big)
        % \nn\\
        &=
            -k_1 e_1^\rmT e_1
            - k_2 e_2^\rmT e_2
            % \nn \\ &\quad
            + e_3^\rmT 
            \xi_3
            % \Big(
            %     \xi_3
            %     +
            %     \kappa_{12}
            %     g_2(- \hat \vartheta_1 + \Theta_1)
            % \Big),
            \label{eq:V3dot_extended}
\end{align}
where 
\begin{align}
    \xi_3 
        &\isdef
        %     2x_2 + k_1 e_{1} + \SG_2 \dot{x}_3\hat{\Theta}_1 + g_2\dot{\hat{\Theta}}_1 
        %     \nn\\& \quad 
        %     + \kappa_{12}
        %     (f_2+g_2{\Theta}_1) 
        %     + k_1k_2x_2,\nn\\
        % & =
            2x_2 + k_1 e_{1} + \SG_2 x_4\hat{\Theta}_1 
            \nn\\ &\quad
            + \gamma_1g_2g_2^\rmT e_2 
            % \nn\\& \quad 
            + \kappa_{12}
            (f_2+g_2{\Theta}_1) 
            + k_1k_2x_2.
        \label{eq:alpha_2}
\end{align}
and 
$\kappa_{12} \isdef k_1-k_2.$
%
% \begin{remark}
%     Note that if
% \begin{align}
%     x_4 
%         &= 
%             -e_3 - 2x_2 - k_1 e_{1} -    \SG_2 \dot x_3\hat{\Theta}_1 \nn\\
%             & \quad - g_2\dot{\hat{\Theta}}_1 - \kappa_{12}(f_2+g_2\Theta_1) - k_1k_2x_2,
% \end{align}
% then, $\dot V_3 <0$.
% \end{remark}
%
% Let $\hat{\vartheta}_1$ be another estimate of $\Theta_1$. 
% Next, define
% \begin{align}
%     x_{4\rmd} 
%         &\isdef 
%             x_4 - 
%             % \Big(
%             %     2x_2 + k_1 e_{1} + \SG_2 \dot x_3\hat{\Theta}_1 + g_2\dot{\hat{\Theta}}_1 
%             %     \nn\\ & \quad
%             %     + \kappa_{12}(f_2+g_2\hat{\vartheta}_1) + k_1k_2x_2
%             % \Big)
%             % \nn\\ & \quad
%             \Big(2x_2 + k_1 e_{1} + \SG_2 \dot x_3\hat{\Theta}_1 + g_2\dot{\hat{\Theta}}_1 
%             \nn\\& \quad 
%             + \kappa_{12}
%             (f_2+g_2{\hat{\vartheta}}_1) 
%             + k_1k_2x_2.\Big) \nn\\
%             & \quad - k_3e_3,\nn\\
%             % \nn\\
%         % & = 
%         %     x_4 - \Big(\kappa_{12} g_2\tilde \vartheta_1 +
%         %     % 2x_2 + k_1 e_{1} 
%         %     % \nn\\& \quad
%         %     % + \SG_2 \dot x_3\hat{\Theta}_1 + g_2\dot{\hat{\Theta}}_1 + \kappa_{12}(f_2+g_2\hat{\vartheta}_1) \nn\\
%         %     % & \quad + k_1k_2x_2
%         %     \alpha
%         %     \Big) 
%         %     - k_3e_3,
%         %     \label{eq:x4d}
%         & =
%             x_4 - \Big(\kappa_{12}g_2\tilde \vartheta_1 + \xi_3\Big)\nn\\
%             & \quad - k_3e_3,
%     \label{eq:x4d_def}
% \end{align}

Let $\hat \vartheta_1$ be an estimate of $\Theta_1,$ and define 
$\tilde \vartheta_1 \isdef \hat \vartheta_1  - \Theta_1.$
Next, define
\begin{align}
    x_{4\rmd} 
        &\isdef
            x_4 - \Big(\kappa_{12}g_2
            \tilde \vartheta_1
            + \xi_3\Big)
            % \nn\\& \quad
            - k_3e_3,
    \label{eq:x4d_def}
\end{align}
where $k_3 > 0,$ and
$\hat{\vartheta}_1$ is given by the adaptation law \eqref{eq:vartheta_1_adaptation}.
% where $\hat{\vartheta}_1 \in \BBR,$
% Let $\hat{\vartheta}_1$ be another estimate of $\Theta_1$.
% which implies that
% \begin{align}
%     \xi_3 
%         &\isdef 2x_2 + k_1 e_{1} + \SG_2 \dot x_3\hat{\Theta}_1 + g_2\dot{\hat{\Theta}}_1 \nn\\
%         & \quad + \kappa_{12}(f_2+g_2\hat{\vartheta}_1) + k_1k_2x_2.
%         \label{eq:V3dot_aux_term}
% \end{align}
% Using \eqref{eq:V3dot_aux_term} in \eqref{eq:x4d} and rearranging, we obtain
Rearranging \eqref{eq:x4d_def} yields
\begin{align}
     \xi_3 
        &= 
            e_4 - \kappa_{12} 
            g_2\tilde \vartheta_1
        % \nn\\& \quad
        % \nn\\ &\quad 
        - k_3e_3 \, \in \mathbb{R}^2,
        \label{eq:V3dot_RHS_term}
\end{align}
where 
\begin{align}
    e_4 \isdef x_4 - x_{4\rmd}.
\end{align}
Substituting \eqref{eq:V3dot_RHS_term} in \eqref{eq:V3dot_extended} yields
\begin{align}
    \dot V_3
        &=
            -k_1e_1^\rmT e_1 - k_2e_2^\rmT e_2 
            - k_3e_3^\rmT e_3 
            + e_3^\rmT e_4 \nn\\
        & \quad - e_3^\rmT \kappa_{12} g_2
        \tilde \vartheta_1
        % (\Theta_1 - \hat{\vartheta}_1)
        .
\end{align}

\subsubsection{$\vartheta_1$ Adaptation}

Next, consider the function
\begin{align}
    \SV_3
        &\isdef
            V_3 + \dfrac{1}{2}\gamma_2^{-1}
            \tilde \vartheta_1 ^\rmT 
            \tilde \vartheta_1 ,
            \label{eq:SV3}
\end{align}
where $\gamma_2>0.$ Differentiating \eqref{eq:SV3} yields
\begin{align}
    \dot{\SV}_3 
        % &= 
        %     -k_1e_1^\rmT e_1 - k_2e_2^\rmT e_2 
        %     - k_3e_3^\rmT e_3 
        %     + e_3^\rmT e_4 
        %     \nn\\ & \quad
        %     - e_3^\rmT \kappa_{12} g_2
        %     \tilde \vartheta_1 
        %     % \nn\\ & \quad
        %     +\gamma_2^{-1} \tilde \vartheta_1 ^\rmT \dot{\hat{\vartheta}}_1 \nn\\
        &=
            -k_1e_1^\rmT e_1 - k_2e_2^\rmT e_2 
            - k_3e_3^\rmT e_3 
            + e_3^\rmT e_4 
            \nn\\ & \quad 
            - \Big(\kappa_{12}g_2^\rmT e_3 -\gamma_2^{-1}\dot{\hat{\vartheta}}_1\Big) \tilde \vartheta_1 .
\end{align}
Letting
\begin{align}
    \dot{\hat{\vartheta}}_1 = \gamma_2\kappa_{12}g_2(x_3)^\rmT e_3,
    \label{eq:vartheta_1_adaptation}
\end{align}
yields
\begin{align}
    \dot{\SV}_3
        &=
            -k_1e_1^\rmT e_1 
            -k_2e_2^\rmT e_2 
            % \nn\\ & \quad
            -k_3e_3^\rmT e_3
            % \nn\\ & \quad
            +e_3^\rmT e_4.
\end{align}

\subsubsection{$e_4$ Stabilization}
Next, consider the function
\begin{align}
    V_4
        &\isdef
            \SV_3 + \dfrac{1}{2}e_4^\rmT e_4.
        \label{eq:V4}
\end{align}
Differentiating \eqref{eq:V4} and using \eqref{eq:x4dot} yields
\begin{align}
 \dot V_4 
        % &= 
        %     -k_1e_1^\rmT e_1 - k_2e_2^\rmT e_2 - k_3e_3^\rmT e_3 
        %     % \nn\\& \quad 
        %     + e_4^\rmT (\dot x_4 - \dot x_{\rm 4d})
        % \nn \\
        % &= 
        %     -k_1e_1^\rmT e_1 - k_2e_2^\rmT e_2 - k_3e_3^\rmT e_3 
        %     \nn\\& \quad 
        %     + e_4^\rmT 
        %     \dfrac{\rmd}{\rmd t} 
        %     (\kappa_{12}g_2\tilde \vartheta_1 + \xi_3 
        %     + k_3e_3)
        % \nn \\
        &= 
            -k_1e_1^\rmT e_1 - k_2e_2^\rmT e_2 - k_3e_3^\rmT e_3 
            % \nn\\& \quad 
            + e_4^\rmT \xi_4
            ,
 \end{align}
 where 
 \begin{align}
     \xi_4 
        &\isdef 
            \Bigg(x_1 + 3f_2 + g_2\hat{\Theta}_{1} + 2k_1x_2
            \nn\\& \quad 
            + k_2\left(x_2 + k_1 e_1\right) + 2g_2\Theta_1
            \nn\\& \quad
            + \dfrac{d}{dt}\left[\dfrac{\partial g(x_3)}{\partial x_3}\right]x_4\hat{\Theta}_{1} 
            \nn\\& \quad 
            + \SG_2 \text{diag}(1,\Theta_2)\hat{\Theta}_{1}u + 2\SG_2 x_4\dot{\hat{\Theta}}_{1}
            \nn\\& \quad
            + \gamma_1g_2\Bigg[x_4^\rmT \SG_2 ^\rmT \left(x_2 + k_1 e_1\right) 
            \nn\\& \quad
            + g_2^\rmT \left(f_2 + g_2\Theta_1 + k_1x_2\right)\Bigg] \nn\\& \quad
            + K_{12} \SG_2 x_4\hat{\vartheta}_{1} 
            \nn\\& \quad 
            + K_{12} g_2\dot{\hat{\vartheta}}_{1} + k_1k_2\left(f_2 + g_2\Theta_1\right) \nn\\& \quad
            + k_3\Bigg(x_2 + \SG_2 x_4\hat{\Theta}_{1} + g_2\dot{\hat{\Theta}}_{1} 
            \nn\\& \quad
            + K_{12} \left(f_2 + g_2\Theta_1\right) + k_1k_2x_2\Bigg)\Bigg)
 \end{align}
 where $K_{12} \isdef k_1 + k_2.$

    Note that if 
    \begin{align}
        \xi_4 = -k_4e_4,
    \end{align}
    where $k_4>0$.
    Then,
    \begin{align}
        \dot V_4
            &=
              -k_1e_1^\rmT e_1 - k_2e_2^\rmT e_2 - k_3e_3^\rmT e_3 
                % \nn\\& \quad 
              -k_4 e_4^\rmT e_4 < 0  
            .
    \end{align}

Let $\hat \varphi_1$ and $\hat \Theta_2$ be an estimate of $\Theta_1$ and $\Theta_2,$ respectively.  
Define 
$\tilde \varphi_1 \isdef \hat \varphi_1  - \Theta_1$ and 
$\tilde \Theta_2 \isdef \hat{\Theta}_2 - \Theta_2.$
Then, letting
\begin{align}
    u
        =& 
            -\left(
            \SG_2 
            \text{diag}(1,\hat{\Theta}_2)\hat{\Theta}_{1}\right)^{-1}
           \Bigg(k_4e_4 + e_3 
           \nn\\ & \quad
           + \left(k_3(k_1+k_2) + k_1k_2 + 2 + \gamma_1g_2~g_2^{\rm T}\right)\left(f_2 + g_2\hat{\varphi}_{1}\right) 
           \nn\\ & \quad
           + \left(k_3(k_1k_2 + 1) + k_1 + \gamma_1g_2~g_2^{\rm T}\right)x_2
           \nn\\ & \quad
           + \left(\SG_2 x_4 -k_3g_2\right)\dot{\hat{\Theta}}_{1} 
           % \nn\\ & \quad
           + \left(  \frac{\rmd} {\rm dt} {\SG_2} + k_3\SG_2 \right)x_4\hat{\Theta}_{1} 
           \nn\\ & \quad
           + K_{12} \SG_2 x_4\hat{\vartheta}_{1}
           + K_{12} g_2\dot {\hat{\vartheta}}_{1} 
            \nn\\& \quad
            +2\gamma_1\SG_2 x_4g_2^{\rm T}(x_2 + k_1 e_1)
            \Bigg)
            ,
            \label{eq:u}
\end{align}
yields
\begin{align}
    \dot V_4 
        % &= 
        %     \dot \SV_3 + e_4^\rmT \dot e_4
        % \nn\\
        % &=
        %     -k_1e_1^\rmT e_1 - k_2e_2^\rmT e_2 
        %     -k_3e_3^\rmT e_3 
        %     + e_4^\rmT  ( e_3 + \dot e_4)
        % \nn \\
        % &=
        %     - k_1e_1^\rmT e_1 
        %     - k_2e_2^\rmT e_2 
        %     - k_3e_3^\rmT e_3
        %     - k_4e_4^\rmT e_4
        %     \nn\\ & \quad % <----- this is easy to comment out if needed. 
        %     + e_4^\rmT \Big[2 + \gamma_1g_2g_2^\rmT  + k_1k_2 
        %     % \nn\\ & \quad 
        %     + k_3K_{12} \Big]g_2\left(\Theta_1 - \hat{\varphi}_{1}\right)\nn\\& \quad           + e_4^\rmT \Big[\SG_2 g_4(x)\text{diag}(\Theta_2,\Theta_2)\hat{\Theta}_{1}u\Big]\left(\Theta_2 - \hat{\Theta}_{2}\right).\nn\\
        % &=
        %     -k_1x_1^{\rm T}x_1 - k_2e_2^{\rm T}e_2 - k_3e_3^{\rm T}e_3 - k_3e_4^{\rm T}e_4\nn\\
        %     & \quad + e_4^{\rm T}\Bigg[\left(k_3K_{12}  + k_1k_2 + 2+ \gamma_1g_2~g_2^{\rm T}\right)g_2\left(\Theta_1-\hat{\varphi}_1\right)\nn\\
        %     & \quad + \left(\SG_2 \text{diag}(0,\tilde \Theta_2)
        %     \hat{\Theta}_{1}u\right)\Bigg],\nn\\
        &=
            -k_1 e_1^{\rm T}e_1 - k_2e_2^{\rm T}e_2 - k_3e_3^{\rm T}e_3 - k_3e_4^{\rm T}e_4
            + e_4^{\rm T} \cdot
            \nn\\ & \quad
            \Bigg[
            -\left(k_3K_{12}  + k_1k_2 + 2+ \gamma_1g_2~g_2^{\rm T}\right)g_2
            % \left(\Theta_1-\hat{\varphi}_1\right)
            \tilde{\varphi}_1
            \nn\\ & \quad 
            - \hat{\Theta}_{1}\SG_2 
            \matl
                0\\
                u_2
            \matr
             \tilde \Theta_2 \Bigg]
\end{align}

\subsubsection{$\hat{\varphi}_1$ and $\hat{\Theta}_2$ Adaptations}
Next, consider the function
\begin{align}
    \SV_4 
        &= 
            V_4
            +
            \dfrac{1}{2}\gamma_3^{-1}
            % (\Theta_1 - \hat{\varphi}_1)^{\rm T}
            % (\Theta_1 - \hat{\varphi}_1)
            \tilde{\varphi}_1^\rmT
            \tilde{\varphi}_1
            % \nn\\ & \quad 
            + \dfrac{1}{2}\gamma_4^{-1} \tilde \Theta_2 ^{\rm T} \tilde \Theta_2 ,
    \label{eq:SV4}
\end{align}
where $\gamma_3>0$ and $\gamma_4>0$. Differentiating \eqref{eq:SV4} yields
\begin{align}
    &\dot{\SV}_4
        = 
            -k_1 e_1^{\rm T}e_1 - k_2e_2^{\rm T}e_2 - k_3e_3^{\rm T}e_3 - k_3e_4^{\rm T}e_4
            \nn\\ &
            % \quad
            + 
            \left(
                -e_4^{\rm T}\big(k_3K_{12}  + k_1k_2 + 2 + \gamma_1g_2~g_2^{\rm T}\big)g_2 
                % \nn\\ & \quad 
                +\gamma_3^{-1}\dot{\hat{\varphi}}_1^{\rm T}
            \right)
            % \left(\Theta_1-\hat{\varphi}_1\right)
            \tilde{\varphi}_1
            \nn\\ & 
            % \quad 
            + 
            \left(
                -e_4^{\rm T}\hat{\Theta}_{1}\SG_2 
                \matl
                    0\\
                    u_2
                \matr
                + \gamma_4^{-1}\dot{\hat{\Theta}}_2^{\rm T}
            \right)
            \tilde \Theta_2 .
            \label{eq:final_V_dot_expanded}
\end{align}
Letting
\begin{align}
    \dot{\hat{\varphi}}_1 
        &= 
            \gamma_3g_2^{\rm T}\big(k_3K_{12}  + k_1k_2 + 2 + \gamma_1g_2^{\rm T}g_2\big)e_4,
        \label{eq:varphi_1_hat_dot}
    \\
    \dot{\hat{\Theta}}_2
        &=
            \gamma_4\hat{\Theta}_{1}
            \matl
                0 & u_2
            \matr
            \SG_2 ^{\rm T}e_4,
        \label{eq:Theta_2_hat_dot}
\end{align}
yields
\begin{align}
    \dot{\SV}_4
        &= 
            -k_1 e_1^{\rm T}x_1 - k_2e_2^{\rm T}e_2 - k_3e_3^{\rm T}e_3 - k_4e_4^{\rm T}e_4.
        \label{eq:final_V_dot}
\end{align}

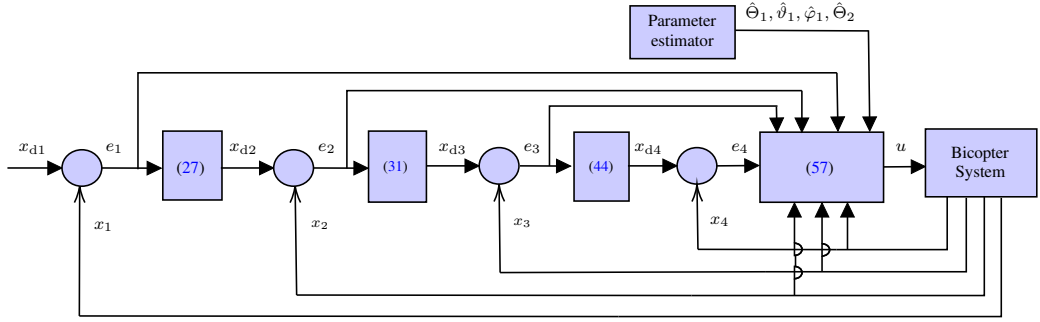
\begin{figure*}[!ht]
    \centering
    \resizebox{1.75\columnwidth}{!}{
    
\tikzset{every picture/.style={line width=0.75pt}} %set default line width to 0.75pt        

\begin{tikzpicture}[x=0.75pt,y=0.75pt,yscale=-1,xscale=1]
%uncomment if require: \path (0,533); %set diagram left start at 0, and has height of 533

%Straight Lines [id:da34647041770057396] 
\draw    (40.95,218.2) -- (72.26,218.2) ;
\draw [shift={(75.26,218.2)}, rotate = 180] [fill={rgb, 255:red, 0; green, 0; blue, 0 }  ][line width=0.08]  [draw opacity=0] (8.93,-4.29) -- (0,0) -- (8.93,4.29) -- cycle    ;
%Shape: Rectangle [id:dp032635824254872325] 
\draw  [fill=blue!20  ,fill opacity=1 ] (138.62,195) -- (175.81,195) -- (175.81,238.33) -- (138.62,238.33) -- cycle ;
%Straight Lines [id:da07104942103763467] 
\draw    (100.26,218.2) -- (134.67,218.32) ;
\draw [shift={(137.67,218.33)}, rotate = 180.2] [fill={rgb, 255:red, 0; green, 0; blue, 0 }  ][line width=0.08]  [draw opacity=0] (8.93,-4.29) -- (0,0) -- (8.93,4.29) -- cycle    ;
%Straight Lines [id:da8971342777114826] 
\draw    (175.09,218.2) -- (205.32,218.2) ;
\draw [shift={(208.32,218.2)}, rotate = 180] [fill={rgb, 255:red, 0; green, 0; blue, 0 }  ][line width=0.08]  [draw opacity=0] (8.93,-4.29) -- (0,0) -- (8.93,4.29) -- cycle    ;
%Shape: Rectangle [id:dp7037990488175827] 
\draw  [fill=blue!20  ,fill opacity=1 ] (267.94,195) -- (304.32,195) -- (304.32,239.67) -- (267.94,239.67) -- cycle ;
%Straight Lines [id:da31783592740818123] 
\draw    (233.32,218.2) -- (264.67,218.32) ;
\draw [shift={(267.67,218.33)}, rotate = 180.22] [fill={rgb, 255:red, 0; green, 0; blue, 0 }  ][line width=0.08]  [draw opacity=0] (8.93,-4.29) -- (0,0) -- (8.93,4.29) -- cycle    ;
%Straight Lines [id:da1954681973345178] 
\draw    (304.32,216.87) -- (334.49,216.87) ;
\draw [shift={(337.49,216.87)}, rotate = 180] [fill={rgb, 255:red, 0; green, 0; blue, 0 }  ][line width=0.08]  [draw opacity=0] (8.93,-4.29) -- (0,0) -- (8.93,4.29) -- cycle    ;
%Shape: Rectangle [id:dp18692737663918368] 
\draw  [fill=blue!20  ,fill opacity=1 ] (397.11,195.67) -- (431.67,195.67) -- (431.67,239) -- (397.11,239) -- cycle ;
%Straight Lines [id:da7740537952488573] 
\draw    (363.14,216.87) -- (392,216.99) ;
\draw [shift={(395,217)}, rotate = 180.24] [fill={rgb, 255:red, 0; green, 0; blue, 0 }  ][line width=0.08]  [draw opacity=0] (8.93,-4.29) -- (0,0) -- (8.93,4.29) -- cycle    ;
%Shape: Ellipse [id:dp5216156749428744] 
\draw  [fill=blue!20  ,fill opacity=1 ] (461.88,216.87) .. controls (461.88,210.54) and (467.48,205.4) .. (474.38,205.4) .. controls (481.29,205.4) and (486.88,210.54) .. (486.88,216.87) .. controls (486.88,223.2) and (481.29,228.33) .. (474.38,228.33) .. controls (467.48,228.33) and (461.88,223.2) .. (461.88,216.87) -- cycle ;
%Straight Lines [id:da05040696974833736] 
\draw    (431,216.87) -- (458.88,216.87) ;
\draw [shift={(461.88,216.87)}, rotate = 180] [fill={rgb, 255:red, 0; green, 0; blue, 0 }  ][line width=0.08]  [draw opacity=0] (8.93,-4.29) -- (0,0) -- (8.93,4.29) -- cycle    ;
%Straight Lines [id:da3845868010282305] 
\draw    (486.88,216.87) -- (510.35,216.67) ;
\draw [shift={(513.35,216.64)}, rotate = 179.51] [fill={rgb, 255:red, 0; green, 0; blue, 0 }  ][line width=0.08]  [draw opacity=0] (8.93,-4.29) -- (0,0) -- (8.93,4.29) -- cycle    ;
%Shape: Rectangle [id:dp48105038410162604] 
\draw  [fill=blue!20  ,fill opacity=1 ] (514,195.67) -- (591.67,195.67) -- (591.67,239.77) -- (514,239.77) -- cycle ;
%Straight Lines [id:da9447773978832872] 
\draw    (591.55,217.53) -- (615.02,217.34) ;
\draw [shift={(618.02,217.31)}, rotate = 179.51] [fill={rgb, 255:red, 0; green, 0; blue, 0 }  ][line width=0.08]  [draw opacity=0] (8.93,-4.29) -- (0,0) -- (8.93,4.29) -- cycle    ;
%Shape: Rectangle [id:dp44877986756688704] 
\draw  [fill=blue!20  ,fill opacity=1 ] (618,193.67) -- (688,193.67) -- (688,236.33) -- (618,236.33) -- cycle ;
%Shape: Right Angle [id:dp6363075123074133] 
\draw   (631.67,236.33) -- (631.67,269.56) -- (562,269.56) ;
%Straight Lines [id:da6982250615344847] 
\draw    (475,269.39) -- (562,269.56) ;
%Straight Lines [id:da13211065667524857] 
\draw    (475,269.39) -- (474.41,230.33) ;
\draw [shift={(474.38,228.33)}, rotate = 89.14] [color={rgb, 255:red, 0; green, 0; blue, 0 }  ][line width=0.75]    (10.93,-3.29) .. controls (6.95,-1.4) and (3.31,-0.3) .. (0,0) .. controls (3.31,0.3) and (6.95,1.4) .. (10.93,3.29)   ;
%Shape: Right Angle [id:dp4528684194886403] 
\draw   (643.67,237) -- (643.67,283.48) -- (517.72,283.48) ;
%Straight Lines [id:da18151225315989028] 
\draw    (350.64,284) -- (517.72,283.48) ;
%Straight Lines [id:da23479835759491885] 
\draw    (350.64,284) -- (350.01,230.33) ;
\draw [shift={(349.99,228.33)}, rotate = 89.33] [color={rgb, 255:red, 0; green, 0; blue, 0 }  ][line width=0.75]    (10.93,-3.29) .. controls (6.95,-1.4) and (3.31,-0.3) .. (0,0) .. controls (3.31,0.3) and (6.95,1.4) .. (10.93,3.29)   ;
%Shape: Right Angle [id:dp5993722597310727] 
\draw   (655,236.33) -- (655,298.44) -- (470.95,298.44) ;
%Straight Lines [id:da8528347499833764] 
\draw    (222.74,298.13) -- (470.95,298.44) ;
%Straight Lines [id:da19982903451598788] 
\draw    (222.74,298.13) -- (221.8,231.03) ;
\draw [shift={(221.77,229.03)}, rotate = 89.19] [color={rgb, 255:red, 0; green, 0; blue, 0 }  ][line width=0.75]    (10.93,-3.29) .. controls (6.95,-1.4) and (3.31,-0.3) .. (0,0) .. controls (3.31,0.3) and (6.95,1.4) .. (10.93,3.29)   ;
%Shape: Right Angle [id:dp9918252918852783] 
\draw   (665.67,237) -- (665.67,311.67) -- (425.02,311.67) ;
%Straight Lines [id:da7866626480403656] 
\draw    (86.15,311.29) -- (425.02,311.67) ;
%Straight Lines [id:da905113795421882] 
\draw    (86.15,311.29) -- (86.36,231.03) ;
\draw [shift={(86.37,229.03)}, rotate = 90.15] [color={rgb, 255:red, 0; green, 0; blue, 0 }  ][line width=0.75]    (10.93,-3.29) .. controls (6.95,-1.4) and (3.31,-0.3) .. (0,0) .. controls (3.31,0.3) and (6.95,1.4) .. (10.93,3.29)   ;
%Straight Lines [id:da11634406093381222] 
\draw    (569,269.67) -- (569,243.14) ;
\draw [shift={(569,240.14)}, rotate = 90] [fill={rgb, 255:red, 0; green, 0; blue, 0 }  ][line width=0.08]  [draw opacity=0] (8.93,-4.29) -- (0,0) -- (8.93,4.29) -- cycle    ;
%Straight Lines [id:da2127563684876983] 
\draw    (553,273.56) -- (553,283.65) ;
%Straight Lines [id:da6601279000877704] 
\draw    (535.67,265.52) -- (535.67,243.14) ;
\draw [shift={(535.67,240.14)}, rotate = 90] [fill={rgb, 255:red, 0; green, 0; blue, 0 }  ][line width=0.08]  [draw opacity=0] (8.93,-4.29) -- (0,0) -- (8.93,4.29) -- cycle    ;
%Shape: Arc [id:dp8273574027765191] 
\draw  [draw opacity=0] (535.67,265.52) .. controls (538.3,265.52) and (540.43,267.25) .. (540.43,269.37) .. controls (540.43,271.49) and (538.3,273.22) .. (535.66,273.22) -- (535.67,269.37) -- cycle ; \draw   (535.67,265.52) .. controls (538.3,265.52) and (540.43,267.25) .. (540.43,269.37) .. controls (540.43,271.49) and (538.3,273.22) .. (535.66,273.22) ;  
%Straight Lines [id:da43864663290303274] 
\draw    (553,265.87) -- (553,243.14) ;
\draw [shift={(553,240.14)}, rotate = 90] [fill={rgb, 255:red, 0; green, 0; blue, 0 }  ][line width=0.08]  [draw opacity=0] (8.93,-4.29) -- (0,0) -- (8.93,4.29) -- cycle    ;
%Shape: Arc [id:dp1969253780109277] 
\draw  [draw opacity=0] (553,265.87) .. controls (555.63,265.87) and (557.76,267.59) .. (557.76,269.72) .. controls (557.76,271.84) and (555.63,273.56) .. (553,273.56) -- (553,269.72) -- cycle ; \draw   (553,265.87) .. controls (555.63,265.87) and (557.76,267.59) .. (557.76,269.72) .. controls (557.76,271.84) and (555.63,273.56) .. (553,273.56) ;  
%Straight Lines [id:da8809457572465964] 
\draw    (535.66,287.83) -- (535.67,298.27) ;
%Shape: Arc [id:dp7510513888653443] 
\draw  [draw opacity=0] (535.67,280.14) .. controls (538.3,280.14) and (540.43,281.86) .. (540.43,283.98) .. controls (540.43,286.11) and (538.3,287.83) .. (535.66,287.83) -- (535.67,283.98) -- cycle ; \draw   (535.67,280.14) .. controls (538.3,280.14) and (540.43,281.86) .. (540.43,283.98) .. controls (540.43,286.11) and (538.3,287.83) .. (535.66,287.83) ;  
%Straight Lines [id:da17629877282207818] 
\draw    (536.33,272.87) -- (536.33,280.17) ;
%Shape: Right Angle [id:dp9494118862368865] 
\draw   (123,219) -- (123,158.33) -- (562.33,158.33) ;
%Straight Lines [id:da7561004525220583] 
\draw    (562.33,158.33) -- (562.95,192) ;
\draw [shift={(563,195)}, rotate = 268.96] [fill={rgb, 255:red, 0; green, 0; blue, 0 }  ][line width=0.08]  [draw opacity=0] (8.93,-4.29) -- (0,0) -- (8.93,4.29) -- cycle    ;
%Shape: Right Angle [id:dp43228530436765555] 
\draw   (254.33,218.49) -- (254.33,169.55) -- (540.33,169.55) ;
%Straight Lines [id:da11622690222739074] 
\draw    (540.33,169.55) -- (540.33,192.67) ;
\draw [shift={(540.33,195.67)}, rotate = 270] [fill={rgb, 255:red, 0; green, 0; blue, 0 }  ][line width=0.08]  [draw opacity=0] (8.93,-4.29) -- (0,0) -- (8.93,4.29) -- cycle    ;
%Shape: Right Angle [id:dp3271349272881743] 
\draw   (381.67,216.96) -- (381.67,179.24) -- (525,179.24) ;
%Straight Lines [id:da5532156489306723] 
\draw    (524.33,179.67) -- (524.6,192.67) ;
\draw [shift={(524.67,195.67)}, rotate = 268.81] [fill={rgb, 255:red, 0; green, 0; blue, 0 }  ][line width=0.08]  [draw opacity=0] (8.93,-4.29) -- (0,0) -- (8.93,4.29) -- cycle    ;
%Shape: Rectangle [id:dp8244288258252093] 
\draw  [fill=blue!20  ,fill opacity=1 ] (432.33,116.22) -- (498.33,116.22) -- (498.33,149.42) -- (432.33,149.42) -- cycle ;
%Straight Lines [id:da0043297319041750715] parameter estimator 
\draw [line width=0.75]    (498.47,131.58) -- (582.17,131.97) ;
%Straight Lines [id:da08485856332954156] 
\draw [line width=0.75]    (582.17,131.58) -- (582.98,192) ;
\draw [shift={(583,195)}, rotate = 269.54] [fill={rgb, 255:red, 0; green, 0; blue, 0 }  ][line width=0.08]  [draw opacity=0] (8.93,-4.29) -- (0,0) -- (8.93,4.29) -- cycle    ;
%Shape: Ellipse [id:dp4654092819541902] 
\draw  [fill=blue!20  ,fill opacity=1 ] (337.49,216.87) .. controls (337.49,210.54) and (343.08,205.4) .. (349.99,205.4) .. controls (356.89,205.4) and (362.49,210.54) .. (362.49,216.87) .. controls (362.49,223.2) and (356.89,228.33) .. (349.99,228.33) .. controls (343.08,228.33) and (337.49,223.2) .. (337.49,216.87) -- cycle ;
%Shape: Ellipse [id:dp6378199806471161] 
\draw  [fill=blue!20  ,fill opacity=1 ] (208.32,218.2) .. controls (208.32,211.87) and (213.92,206.74) .. (220.82,206.74) .. controls (227.72,206.74) and (233.32,211.87) .. (233.32,218.2) .. controls (233.32,224.53) and (227.72,229.67) .. (220.82,229.67) .. controls (213.92,229.67) and (208.32,224.53) .. (208.32,218.2) -- cycle ;
%Shape: Ellipse [id:dp018945879810865707] 
\draw  [fill=blue!20  ,fill opacity=1 ] (75.26,218.2) .. controls (75.26,211.87) and (80.85,206.74) .. (87.76,206.74) .. controls (94.66,206.74) and (100.26,211.87) .. (100.26,218.2) .. controls (100.26,224.53) and (94.66,229.67) .. (87.76,229.67) .. controls (80.85,229.67) and (75.26,224.53) .. (75.26,218.2) -- cycle ;

% Text Node
\draw (101.34,201.25) node [anchor=north west][inner sep=0.75pt]  [font=\footnotesize] [align=left] {$e_1$};
% Text Node
\draw (46.67,201.25) node [anchor=north west][inner sep=0.75pt]  [font=\footnotesize] [align=left] {$x_{\rm{d}1}$};
% Text Node
\draw (145.79,212.47) node [anchor=north west][inner sep=0.75pt]  [font=\footnotesize] [align=left] {\eqref{eq:x2d}};
% Text Node
\draw (234.87,201.25) node [anchor=north west][inner sep=0.75pt]  [font=\footnotesize] [align=left] {$e_2$};
% Text Node
\draw (179.14,201.25) node [anchor=north west][inner sep=0.75pt]  [font=\footnotesize] [align=left] {$x_{\rm{d}2}$};
% Text Node
\draw (275.7,212.47) node [anchor=north west][inner sep=0.75pt]  [font=\scriptsize] [align=left] {\eqref{eq:x3d}};
% Text Node
\draw (363.95,201.25) node [anchor=north west][inner sep=0.75pt]  [font=\footnotesize] [align=left] {$e_3$};
% Text Node
\draw (310.64,201.25) node [anchor=north west][inner sep=0.75pt]  [font=\footnotesize] [align=left] {$x_{\rm{d}3}$};
% Text Node
\draw (404.96,212.47) node [anchor=north west][inner sep=0.75pt]  [font=\scriptsize] [align=left] {\eqref{eq:x4d_def}};
% Text Node
\draw (494.5,201.25) node [anchor=north west][inner sep=0.75pt]  [font=\footnotesize] [align=left] {$e_4$};
% Text Node
\draw (433.67,201.25) node [anchor=north west][inner sep=0.75pt]  [font=\footnotesize] [align=left] {$x_{\rm{d}4}$};
% Text Node
\draw (597.83,201.25) node [anchor=north west][inner sep=0.75pt]  [font=\footnotesize] [align=left] {$u$};
% Text Node
\draw (629.29,203.14) node [anchor=north west][inner sep=0.75pt]  [font=\footnotesize] [align=left] {\begin{minipage}[lt]{32.66pt}\setlength\topsep{0pt}
\begin{center}
Bicopter\\System
\end{center}

\end{minipage}};
% Text Node
\draw (541.29,212.47) node [anchor=north west][inner sep=0.75pt]  [font=\footnotesize] [align=left] {\eqref{eq:u}};
% Text Node
\draw (93.92,250.58) node [anchor=north west][inner sep=0.75pt]  [font=\footnotesize] [align=left] {$x_1$};
% Text Node
\draw (230.31,249.91) node [anchor=north west][inner sep=0.75pt]  [font=\scriptsize] [align=left] {$x_2$};
% Text Node
\draw (357.5,248.58) node [anchor=north west][inner sep=0.75pt]  [font=\scriptsize] [align=left] {$x_3$};
% Text Node
\draw (481.2,247.25) node [anchor=north west][inner sep=0.75pt]  [font=\footnotesize] [align=left] {$x_4$};
% Text Node
\draw (436.62,120.5) node [anchor=north west][inner sep=0.75pt]  [font=\footnotesize] [align=left] {\begin{minipage}[lt]{40.82pt}\setlength\topsep{0pt}
\begin{center}
Parameter\\estimator
\end{center}

\end{minipage}};
% Text Node
\draw (503.47,112.28) node [anchor=north west][inner sep=0.75pt]  [font=\footnotesize] [align=left] {$
\begingroup % keep the change local
\setlength\arraycolsep{2pt}
% \matl
\hat{\Theta}_1 , \hat{\vartheta}_1 , \hat{\varphi}_1 , \hat{\Theta}_2
% \matr
\endgroup$};

\end{tikzpicture}
}
    \vspace{1pt}
    \caption{Adaptive backstepping control architecture.}
    \label{fig:DEA-BACK_architecture}
\end{figure*}

The controller is thus given by 
\eqref{eq:x2d},
\eqref{eq:x3d},
\eqref{eq:x4d_def}, and 
\eqref{eq:u} with the parameter adaption laws 
\eqref{eq:hat_Theta_1_dot}, 
\eqref{eq:vartheta_1_adaptation},
\eqref{eq:varphi_1_hat_dot}, and
\eqref{eq:Theta_2_hat_dot}.
The block diagram showing the architecture of the adaptive backstepping controller is shown in Figure \ref{fig:DEA-BACK_architecture}.

\subsection{Stability Analysis}
Note that the control \eqref{eq:u} requires
$    \SG_2 
    \text{diag}(1,\hat{\Theta}_2)\hat{\Theta}_{1}
$ to be nonsingular, 
%         \neq
%             0,
% \end{align}
which implies that $\SG_2$ must be nonsingular and $\hat{\Theta}_1, \hat{\Theta}_2 \neq 0$.
If $F \neq 0,$ which is reasonably expected during the system's operation, $\SG_2$ is nonsingular. 
% 
% As long as $x_{3,1} \neq 0$, $\SG_2\neq 0 $, this implies that as long as $F \neq 0$, $u$ is bounded. 
% 
However, the parameter adaptation laws do not ensure that $\hat{\Theta}_1, \hat{\Theta}_2 \neq 0.$
Thus, the global asymptotic stability of the closed-loop system can not be guaranteed.

% However, if 

Next, consider the function 
\begin{align}
    &\SV_4(e_1, e_2, e_3, e_4, \hat \Theta_1, \hat \vartheta_1, \hat \varphi_1, \hat \Theta_2 ) 
        = 
            \dfrac{1}{2}e_1^{\rm T}e_1
            +
            \dfrac{1}{2}e_2^{\rm T}e_2 
            \nn\\ &\quad 
            + \dfrac{1}{2}e_3^{\rm T}e_3
            +
            \dfrac{1}{2}e_4^{\rm T}e_4
            % \nn\\ & \quad
            + \dfrac{1}{2}\gamma_1^{-1}\tilde{\Theta}_1^{\rm T}\tilde{\Theta}_1
            % \nn\\ & \quad
            + \dfrac{1}{2}\gamma_2^{-1}\tilde{\vartheta}_1^{\rm T}\tilde{\vartheta}_1
            \nn\\ & \quad
            + \dfrac{1}{2}\gamma_3^{-1}\tilde{\varphi}_1^{\rm T}\tilde{\varphi}_1
            % \nn\\ & \quad
            + \dfrac{1}{2}\gamma_4^{-1}\tilde{\Theta}_2^{\rm T}\tilde{\Theta}_2,
\end{align}
where $\gamma_1, \gamma_2, \gamma_3, \gamma_4 > 0.$
Let $D = \BBR^2 \times \BBR^2 \times \BBR^2 \times \BBR^2 \times \BBR/\{0 \} \times \BBR \times \BBR \times \BBR/\{0 \}.$
Then, $\SV_4(x) \ge 0$ for all $x \in D.$
% Note that $\SV_4 \ge 0$ for all $e_1, e_2, e_3, e_4, \hat \Theta_1, \hat \vartheta_1, \hat \varphi_1, \hat \Theta_2.$
Furthermore, it follows from \eqref{eq:final_V_dot} that $\dot \SV_4 \leq 0.$
Moreover, $\dot \SV_4(x) = 0$ if and only if $e_1 = e_2 = e_3 = e_4 = 0.$
Let $E$ be the set of all points in $D$ such that $e_1 = e_2 = e_3 = e_4 = 0.$
Note that $E$ is the largest positively invariant set with respect to \eqref{eq:x1dot}-\eqref{eq:x4dot} in $E.$ 
Then, it follows from the Barbashin-Krasovsky-La Salle's invariance principle \cite{khalil2013nonlinear} that $e_1, e_2, e_3, e_4 \to 0.$

\section{Simulations}
\label{sec:simulations}

In this section, we apply the adaptive backstepping controller developed in the previous section to the trajectory tracking problem. 
% 
% This section presents the numerical simulation results obtained by applying the adaptive backstepping controller to the trajectory following problem for a bicopter. 
In particular, we use the adaptive controller to follow an elliptical and a  second-order Hilbert curve based trajectory. 

To simulate the bicopter, we assume that the mass of the bicopter is $1 \ \rm kg$ and its inertia is $0.2$ $\rm kg \cdot m^2.$ 
In the controller, 
we set 
$k_1 = 5$, $k_2=5$, $k_3=4$, and $k_4=4$.
In the adaptation laws, 
we set the adaptation gains 
$\gamma_1=1$, $\gamma_2=0.05$, $\gamma_3=0.05$, and $\gamma_4=0.1$,  
and the initial estimates 
% $\hat{\Theta}_1(0) = \matl 0.5 & 0.5 & 0.5\matr^\rmT $, and 
$\hat{\Theta}_1(0) = 0.5 $,
$\hat{\vartheta}_1(0) = 0.5 $,
$\hat{\varphi}_1(0) = 0.5 $,
and 
$\hat{\Theta}_2(0) = 40$.

\subsection{Elliptical Trajectory}
\begin{figure}[!ht]
    \centering
    \includegraphics[width=\linewidth]{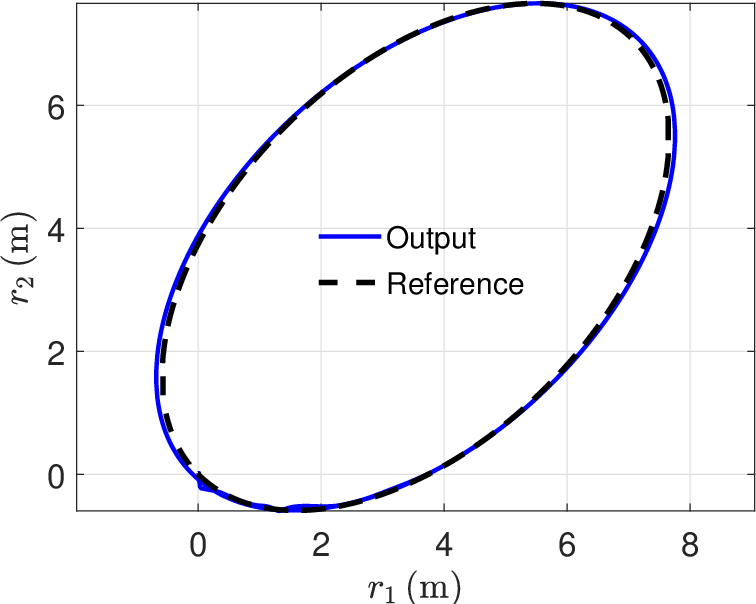}
    \caption{\textbf{Elliptical trajectory}. Tracking response of the bicopter with the adaptive backstepping controller. 
    Note that the output trajectory is shown in solid blue, and the reference trajectory is in dashed black.}
    \label{fig:ACC_DEABACK_Bicopter_Trajectory_Elliptical}
\end{figure}

\begin{figure}[!ht]
    \centering
    \includegraphics[width=\linewidth]{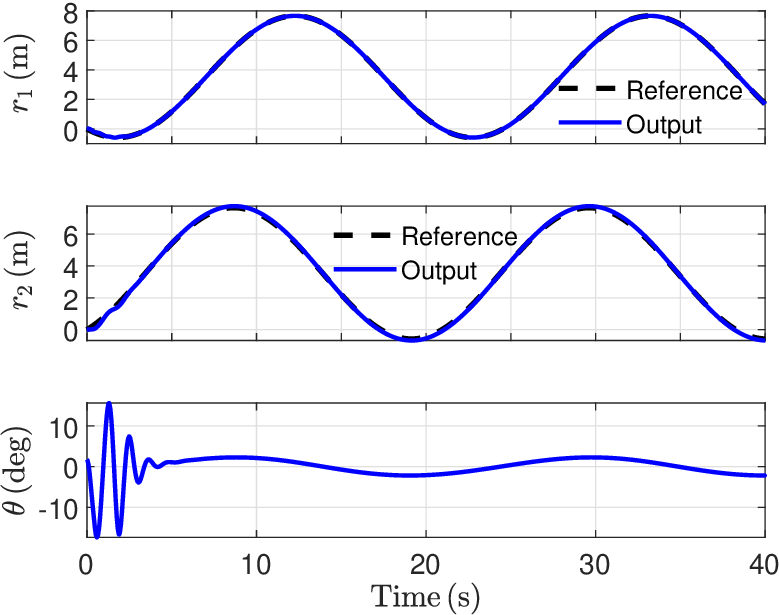}
    \caption{\textbf{Elliptical trajectory}. Position $(r_1, r_2)$ and roll angle $\theta$ response of the bicopter obtained with adaptive backstepping controller \eqref{eq:u}.}
    \label{fig:ACC_DEABACK_Bicopter_states_Elliptical}
\end{figure}

\begin{figure}[!ht]
    \centering
    \includegraphics[width=\linewidth]{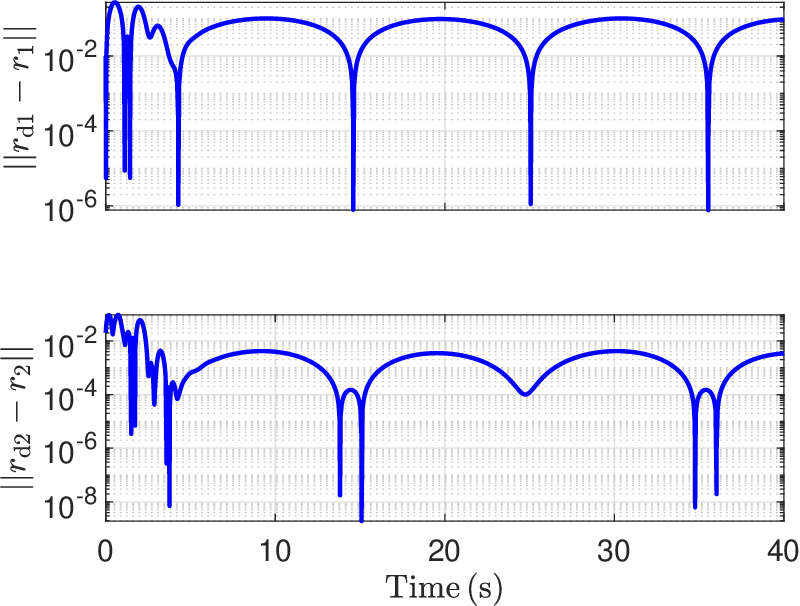}
    \caption{\textbf{Elliptical trajectory}. 
    Position errors obtained with the adaptive backstepping controller  \eqref{eq:u} on a logarithmic scale.}
    \label{fig:ACC_DEABACK_Bicopter_states_errors_Elliptical}
\end{figure}

The bicopter is commanded to follow a elliptical trajectory given by 
\begin{align}
    r_{\rmd1}(t) &= 5 \cos(\phi)-5 \cos(\phi) \cos (\omega t) - 3 \sin(\phi)\sin(\omega t), \nn\\ 
    r_{\rmd2}(t) &= 5\sin(\phi) - 5 \sin(\phi) \cos (\omega t) + 3 \cos(\phi)\sin(\omega t), \nn
\end{align}
where $\phi=45~\rm{deg}$ and $\omega = 0.3 \ \rm rad/s^{-1}. $
Figure \ref{fig:ACC_DEABACK_Bicopter_Trajectory_Elliptical} shows the trajectory-tracking response of the bicopter, where the desired trajectory is shown in black dashes and the output trajectory response is shown in blue.
Figure \ref{fig:ACC_DEABACK_Bicopter_states_Elliptical} shows the position $r_1, r_2$ response and the roll angle $\theta$ response of the bicopter with the adaptive backstepping controller. 
Figure \ref{fig:ACC_DEABACK_Bicopter_states_errors_Elliptical} shows the norm of the position errors obtained with the adaptive backstepping controller on a logarithmic scale. 


\begin{figure}[!ht]
    \centering
    \includegraphics[width=\linewidth]{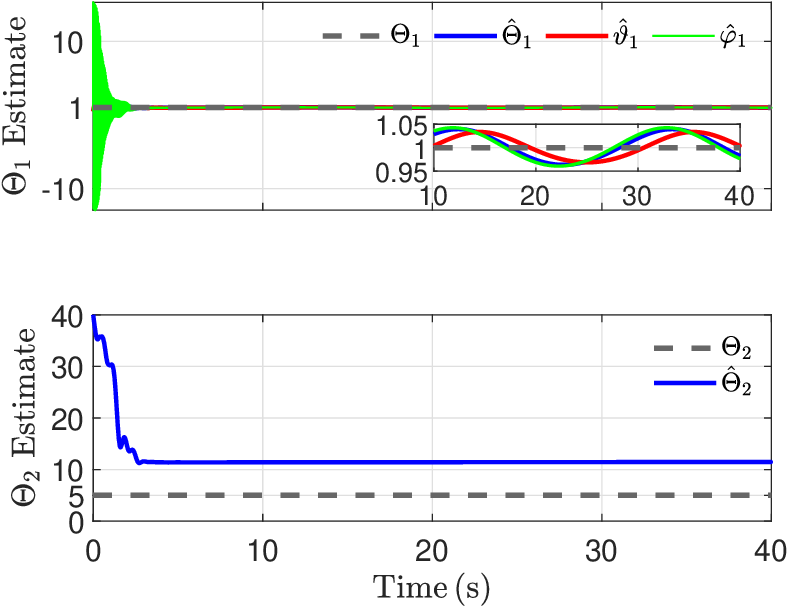}
    \caption{\textbf{Elliptical trajectory}. 
    Estimates of $\Theta_1$ and $\Theta_2$ obtained with adaption laws 
    \eqref{eq:hat_Theta_1_dot}, 
    \eqref{eq:vartheta_1_adaptation},
    \eqref{eq:varphi_1_hat_dot}, and
    \eqref{eq:Theta_2_hat_dot}.
    % The true values parameters used in the model are in dashed black.
    }
    \label{fig:ACC_DEABACK_Bicopter_estimations_Elliptical}
\end{figure}

\begin{figure}[!ht]
    \centering
    \includegraphics[width=\linewidth]{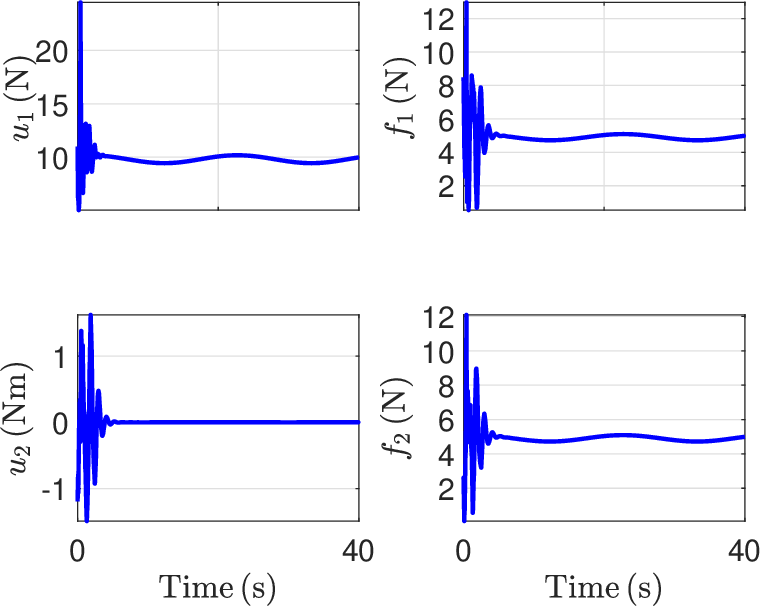}
    \caption{\textbf{Elliptical trajectory}. 
    Control $u$ and the corresponding forces $f_1$ and $f_2$ obtained with adaptive backstepping controller \eqref{eq:u}.}
    \label{fig:ACC_DEABACK_Bicopter_control_Elliptical}
\end{figure}

Figure \ref{fig:ACC_DEABACK_Bicopter_estimations_Elliptical} shows the estimates $\hat \Theta_1,\hat \vartheta_1,$ and $ \hat \varphi_1$ of $\Theta_1$ and the estimate $\hat \Theta_2$ of $\Theta_2.$
Note that the parameter estimates do not converge to their actual values.
However, the non-convergence of the estimates is not due to persistency-related issues. 
In the adaptive controller design, since the parameter adaptation laws are chosen to cancel undesirable factors and not to estimate the parameters, the estimates do not necessarily need to converge.
Finally, Figure \ref{fig:ACC_DEABACK_Bicopter_control_Elliptical} shows the control $u$ generated by the adaptive backstepping controller \eqref{eq:u} and the corresponding forces $f_1$ and $f_2.$
Note that the forces $f_1, f_2$ are computed using \eqref{control_signals_1},and \eqref{control_signals_2}. 


\subsection{Hilbert trajectory}
\begin{figure}[!ht]
    \centering
    \includegraphics[width=\linewidth]{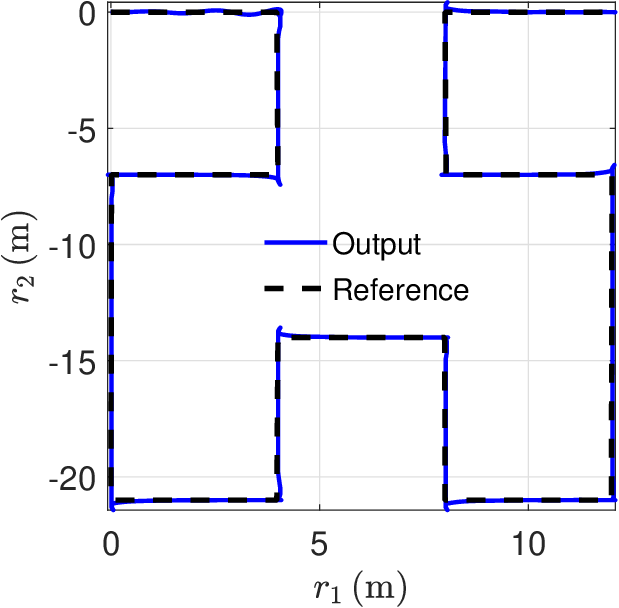}
    \caption{\textbf{Hilbert trajectory}. Tracking response of the bicopter with adaptive backstepping controller. Note that the output trajectory is in solid blue, and the desired trajectory is in dashed black.}
    \label{fig:ACC_DEABACK_Bicopter_Trajectory_Hilbert}
\end{figure}

\begin{figure}[!ht]
    \centering
    \includegraphics[width=\linewidth]{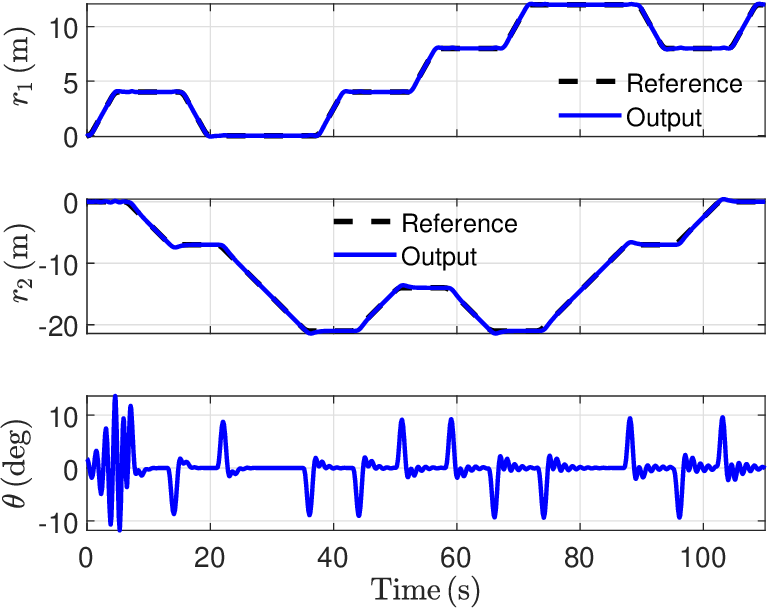}
    \caption{\textbf{Hilbert trajectory}. Positions $(r_1, r_2)$ and roll angle $\theta$ response of the bicopter with adaptive backstepping controller.}
    \label{fig:ACC_DEABACK_Bicopter_states_Hilbert}
\end{figure}

\begin{figure}[!ht]
    \centering
    \includegraphics[width=\linewidth]{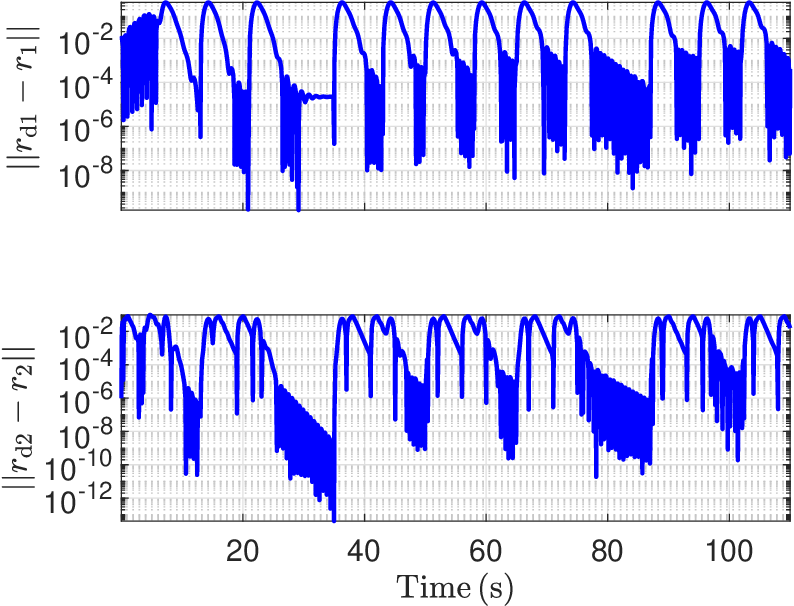}
    \caption{\textbf{Hilbert trajectory}. 
    Position errors with the adaptive backstepping controller on a logarithmic scale.}
    \label{fig:ACC_DEABACK_Bicopter_states_errors_Hilbert}
\end{figure}

Next, the bicopter is commanded to follow a nonsmooth trajectory constructed using a second-order Hilbert curve.
The trajectory is constructed using the algorithm described in Appendix A of \cite{spencer2022adaptive} with
a maximum velocity $v_{\rm max} = 1 \ \rm m/s$ and 
a maximum acceleration $a_{\rm max} = 1 \ \rm m/s^2.$
Figure \ref{fig:ACC_DEABACK_Bicopter_Trajectory_Hilbert} shows the trajectory-tracking response of the bicopter, where the desired trajectory is shown in black dashes, and the output trajectory response is shown in blue. 
Figure \ref{fig:ACC_DEABACK_Bicopter_states_Hilbert} shows the position $r_1$ and $r_2$ response and the roll angle $\theta$ response of the bicopter with the adaptive backstepping controller. 
Figure \ref{fig:ACC_DEABACK_Bicopter_states_errors_Hilbert} shows the norm of the position errors obtained with the adaptive backstepping controller on a logarithmic scale. 


\begin{figure}[!ht]
    \centering
    \includegraphics[width=\linewidth]{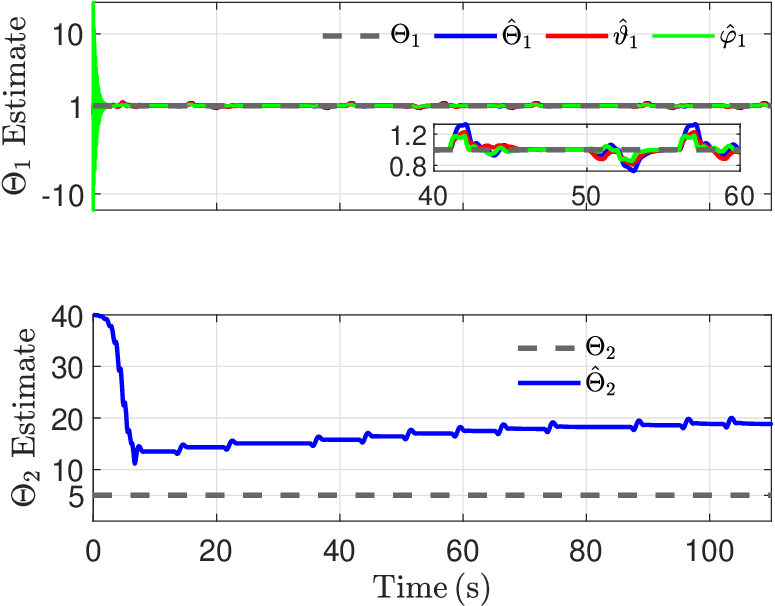}
    \caption{\textbf{Hilbert trajectory}. 
    $\Theta_1$ and $\Theta_2$ Estimations with the adaptive backstepping controller. Note that the parameters used in the model are in dashed black.}
    \label{fig:ACC_DEABACK_Bicopter_estimations_Hilbert}
\end{figure}

\begin{figure}[!ht]
    \centering
    \includegraphics[width=\linewidth]{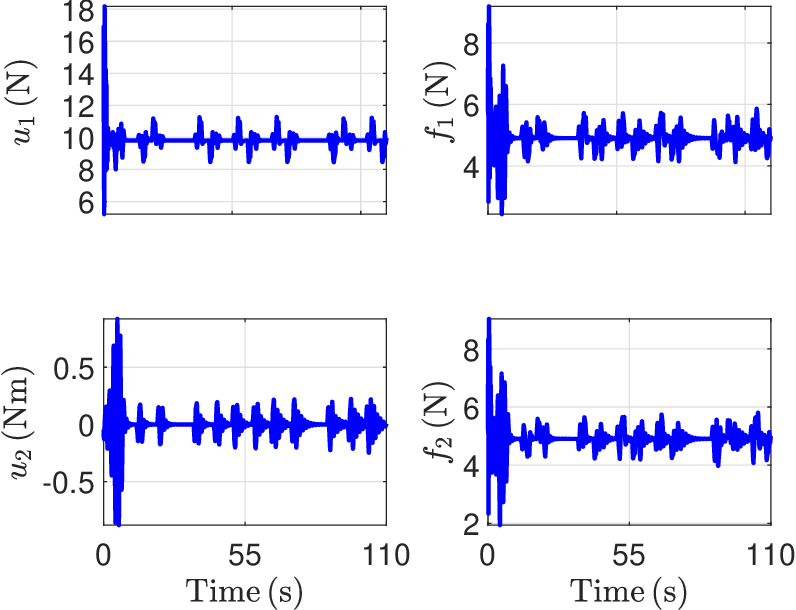}
    \caption{\textbf{Hilbert trajectory}. 
    Control $u$ and the corresponding forces $f_1$ and $f_2$ obtained with adaptive backstepping controller.}
    \label{fig:ACC_DEABACK_Bicopter_control_Hilbert}
\end{figure}

Figure \ref{fig:ACC_DEABACK_Bicopter_estimations_Hilbert} shows the estimates $\hat \Theta_1,\hat \vartheta_1,$ and $\hat \varphi_1$ of $\Theta_1$ and the estimate $\hat \Theta_2$ of $\Theta_2.$
As in the previous case, note that the estimates do not converge to their actual values. 
Finally, Figure \ref{fig:ACC_DEABACK_Bicopter_control_Hilbert} shows the control $u$ generated by the adaptive backstepping controller \eqref{eq:u}, and the corresponding forces $f_1$ and $f_2$. Note that the forces $f_1$ and $f_2$ are computed using \eqref{control_signals_1}, and \eqref{control_signals_2}.


The preceding two examples show that the adaptive backstepping controller stabilizes the bicopter dynamics and successfully tracks the desired trajectory without prior knowledge of the bicopter dynamics.

% bicopter's trajectory converges to the desired trajectory exponentially without prior perfect knowledge of the bicopter dynamics.

\section{Conclusions}
\label{sec:conclusions}
This paper presented an adaptive backstepping-based controller for the stabilization and tracking problem in a bicopter system. 
It is shown that the backstepping process can not be applied to the bicopter dynamics due to the singularity of the input map.
The bicopter dynamics is then dynamically extended to circumvent the singularity problem.
% , but yields a pure feedback system, where one of the virtual control signals appears non-affinely in the system. 
% 
The backstepping process is then used to stabilize the successive states of the extended bicopter system and design parameter adaption laws.
Since the final control law requires the inversion of two of the parameter estimates, the global asymptotic stability of the closed-loop system cannot be guaranteed as the adaptation laws do not enforce a constraint on the estimate values. 
However, as shown in the numerical simulations, the closed-loop system is asymptotically stable with appropriately chosen gains, and the controller yields the desired tracking performance. 

Future extensions will focus on
1) integrating constraints in the parameter adaption law that guarantees global stability of the closed-loop system for all gains chosen in the control law and the parameter adaptation laws,
2) reformulating the control laws as cascaded control, 
and finally,
3) extending the adaptive backstepping control presented in this paper to design a stabilizing and tracking controller for a quadcopter. 

% The controller was constructed by extending the bicopter dynamics to allow expressing the system in pure-feedback triangular form and designing an adaptive backstepping controller to perform tracking of both smooth and non-smooth trajectories.
% \clearpage

\printbibliography
\end{document}